\documentclass[11pt]{amsart}
\usepackage[left=3cm,top=3cm,right=3cm,bottom=3cm]{geometry}
\usepackage{amsfonts}
\usepackage{enumerate,subfigure,color}
\pdfoutput=1
\usepackage{amssymb,latexsym,amsmath}
\usepackage[pdftex]{graphicx}
\usepackage{framed}

\newcommand{\dz}{\partial_z}

\newcommand{\ra}{\rightarrow}

\newcommand{\T}{{\mathbf{t}}}

\newcommand{\R}{{\mathbb R}}

\newcommand{\C}{{\mathbb C}}

\newcommand{\norm}[1]{\left\|#1\right\|}

\DeclareMathOperator{\p}{\partial}

\DeclareMathOperator{\by}{\times}
\DeclareMathOperator{\bndry}{\partial\Omega}

\def\bfR{\mathbb{R}}
\def\bfu{\mbox{\boldmath${u}$}}
\def\bfp{\mbox{\boldmath${p}$}}
\def\bff{\mbox{\boldmath${f}$}}

\def\bfS{\mbox{\boldmath$S$}}
\def\bfT{\mbox{\boldmath$T$}}

\title[Hybrid segmentation and D-bar method for EIT]{A Hybrid Segmentation and D-bar Method for Electrical Impedance Tomography}
\author{S.~J. Hamilton}\thanks{Hamilton: Department of Mathematics, Statistics, and Computer Science; Marquette University, Milwaukee, Wisconsin USA}

\author{J.~M. Reyes}\thanks{Reyes: School of Computer Science \& Informatics; Cardiff University, Cardiff, United Kingdom}
\author{S. Siltanen}\thanks{Siltanen: Department of Mathematics \& Statistics; University of Helsinki, Helsinki, Finland}
\author{X. Zhang}\thanks{Zhang: Department of Mathematics, MOE-LSC, and Institute of Natural Sciences; Shanghai Jiao Tong University, Shanghai, China}
\date{June 11, 2015}
\email{sarah.hamilton@marquette.edu}
\email{reyes.juanmanuel@gmail.com} \email{samuli.siltanen@iki.fi}
\email{xqzhang@sjtu.edu.cn}

\begin{document}
\begin{abstract}
The Regularized D-bar method for Electrical Im\-pe\-dan\-ce
Tomography provides a rigorous mathematical approach for solving
the full nonlinear inverse problem directly, i.e. without iterations. It is based on a low-pass filtering in the (nonlinear) frequency domain.  However, the resulting D-bar reconstructions are inherently smoothed leading to a loss of edge distinction.  In this paper, a novel approach that combines the rigor of the D-bar approach with the edge-preserving nature of Total Variation regularization is presented.  The method also includes a data-driven contrast adjustment technique guided by the key functions (\emph{CGO solutions}) of the D-bar method.  The new \emph{TV-Enhanced D-bar Method} produces reconstructions with sharper edges and improved
contrast while still solving the full nonlinear problem. This is achieved by
using the TV-induced edges to increase the truncation radius of the scattering data in the nonlinear frequency domain thereby increasing the radius of the low pass filter. The algorithm is tested on numerically simulated noisy EIT data and demonstrates significant improvements in edge preservation and
contrast which can be highly valuable for absolute EIT imaging.
\end{abstract}
\maketitle
%
%
\section{Introduction}

\noindent
In \emph{Electrical Impedance Tomography} (EIT) a conductive body is probed with harmless electrical currents fed into the body through electrodes at the surface, and the resulting voltages are measured at the electrodes. The goal is to recover the electrical conductivity distribution inside the body from these surface electrical measurements. EIT is useful in medical imaging, as different tissues have different conductivities, and it allows harmless and painless monitoring of patients even over long periods of time. Another application area of EIT is non-destructive testing.
See \cite{Cheney1999, Mueller2012} for reviews of EIT and its uses.

The image reconstruction task of EIT is a nonlinear and severely ill-posed inverse problem. Therefore, EIT algorithms need to be {\em regularized} to overcome the extreme sensitivity to modeling errors and measurement noise.
Among EIT algorithms, the so-called {\em D-bar method} stands out due to its unique capability of dividing the measurement information neatly into stable and unstable parts in a (nonlinear) frequency domain.
See Figure \ref{fig:scheme}.

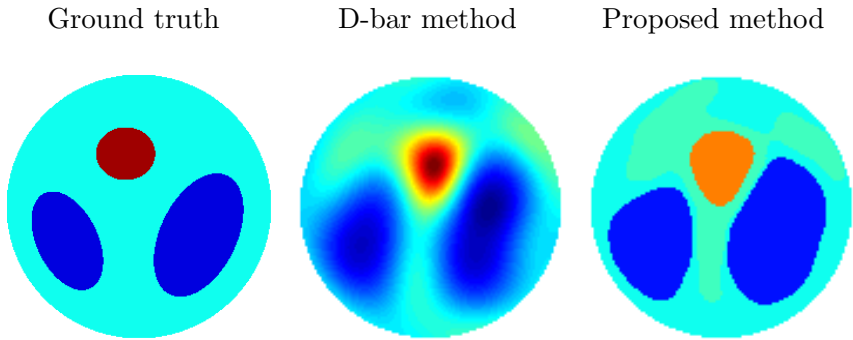
\begin{figure}[t!]
\centering
\begin{picture}(330,130)
\put(0,2){\includegraphics[width=100pt]{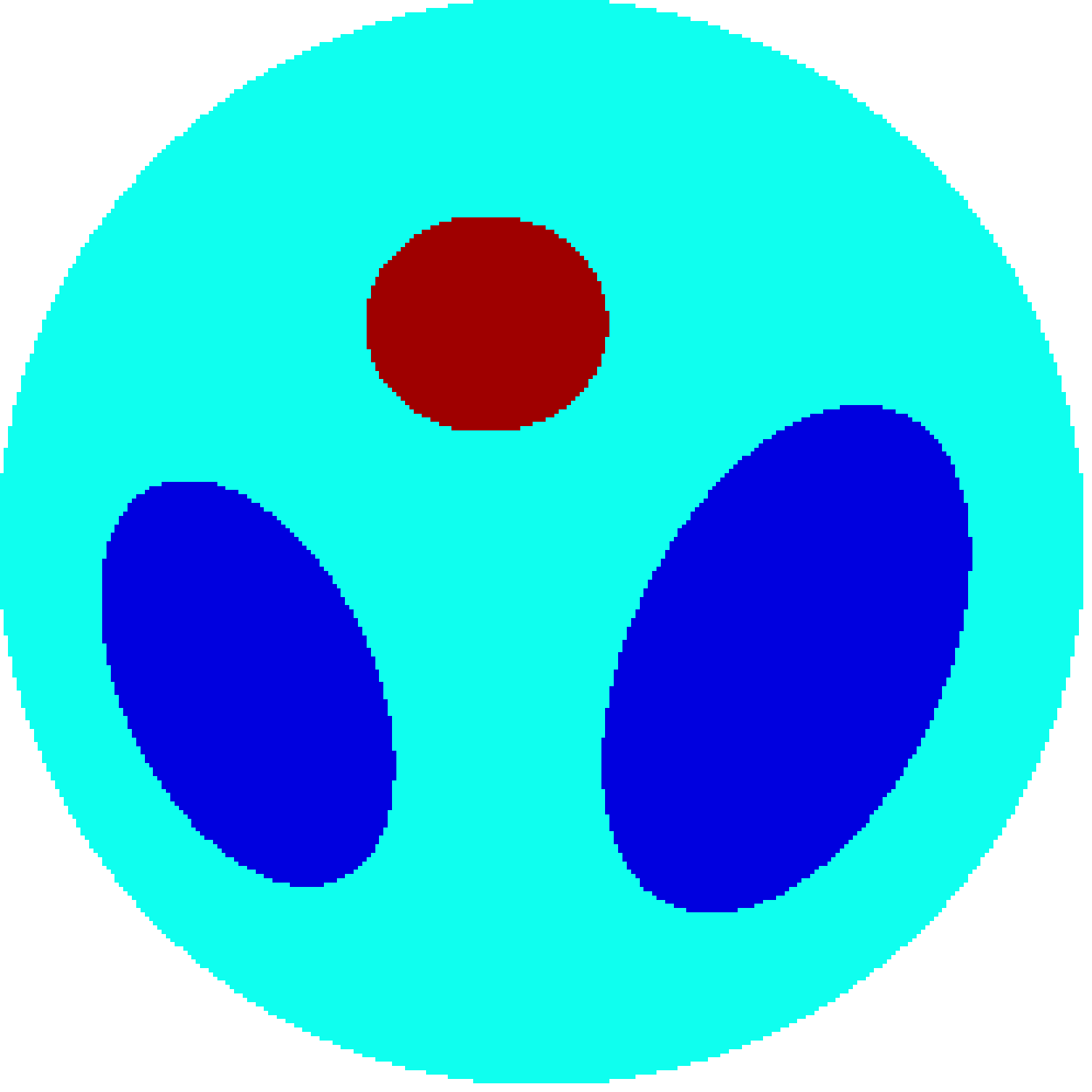}}
\put(110,2){\includegraphics[width=100pt]{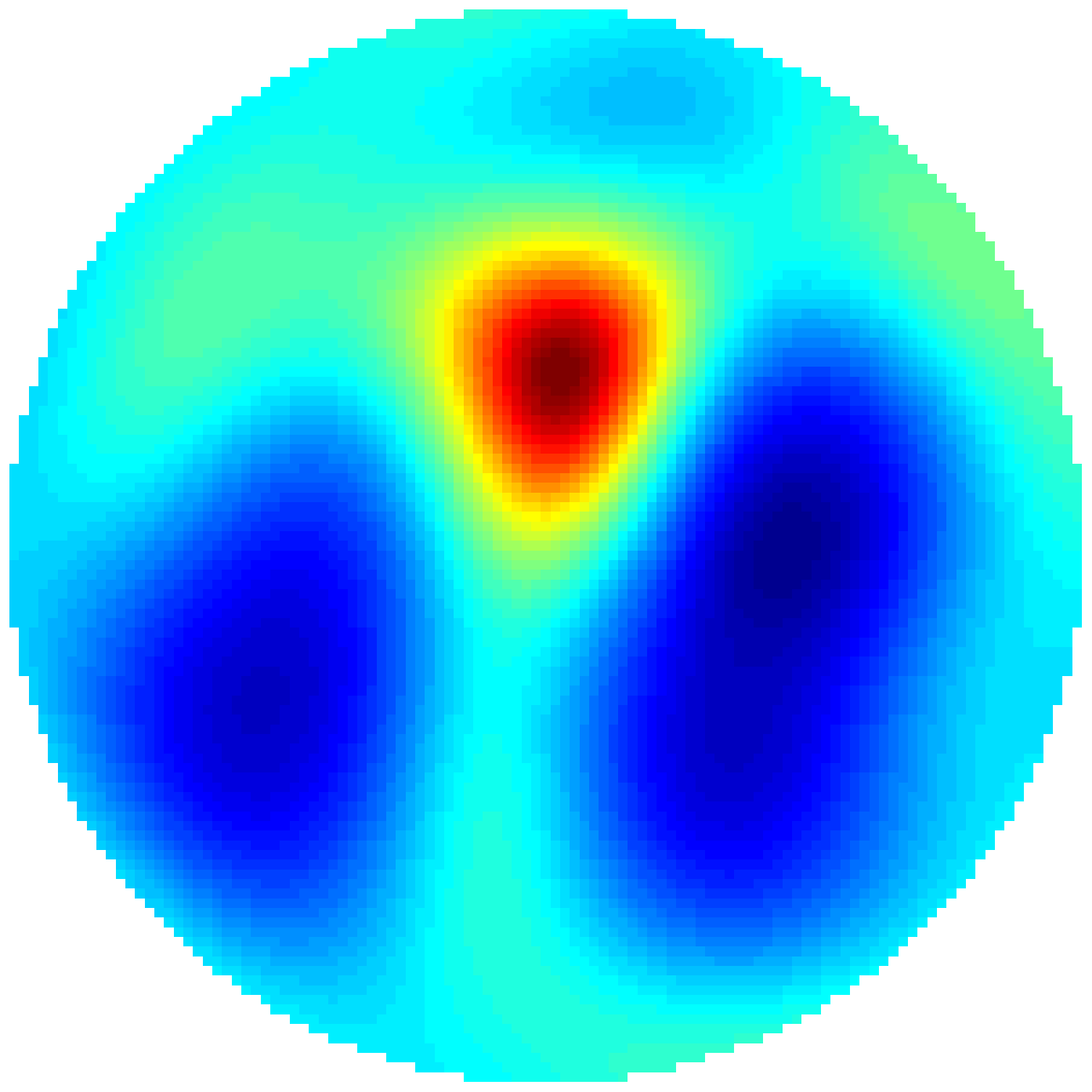}}
\put(220,2){\includegraphics[width=100pt]{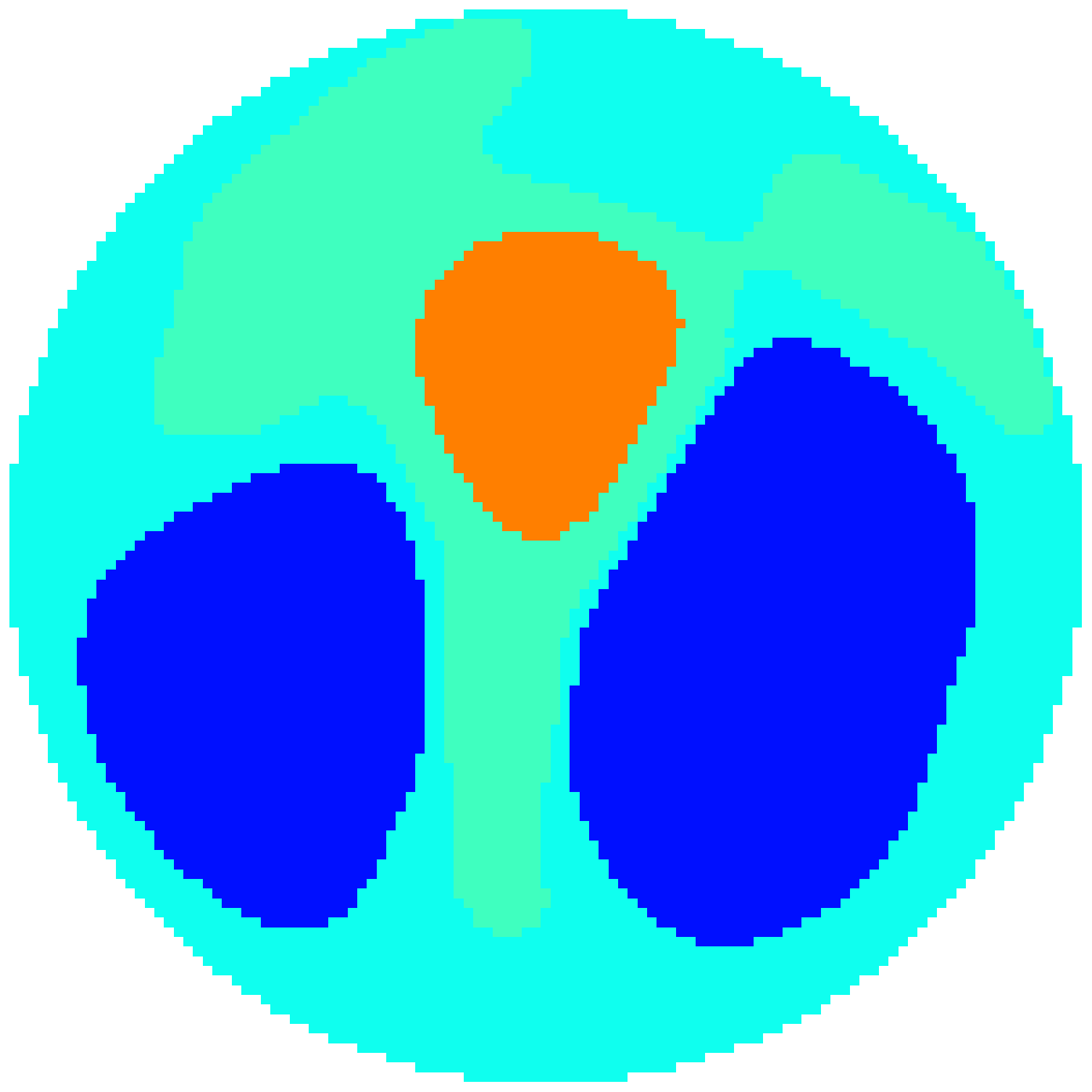}}

\put(15,120){Ground truth}
\put(125,120){D-bar method}
\put(225,120){Proposed method}
\end{picture}
\caption{\label{fig:hammer}Left: simulated ``heart-and-lungs'' phantom conductivity. Middle: D-bar reconstruction based on nonlinear low-pass filtering with 0.75\% relative noise added to EIT voltage data (see Figure \ref{fig:scheme}). Right: reconstruction with the proposed hybrid method from the same noisy EIT data.}
\end{figure}

\begin{figure}[t!]
\begin{center}
\begin{picture}(320,240)
\put(-30,2){\includegraphics[height=2.5cm]{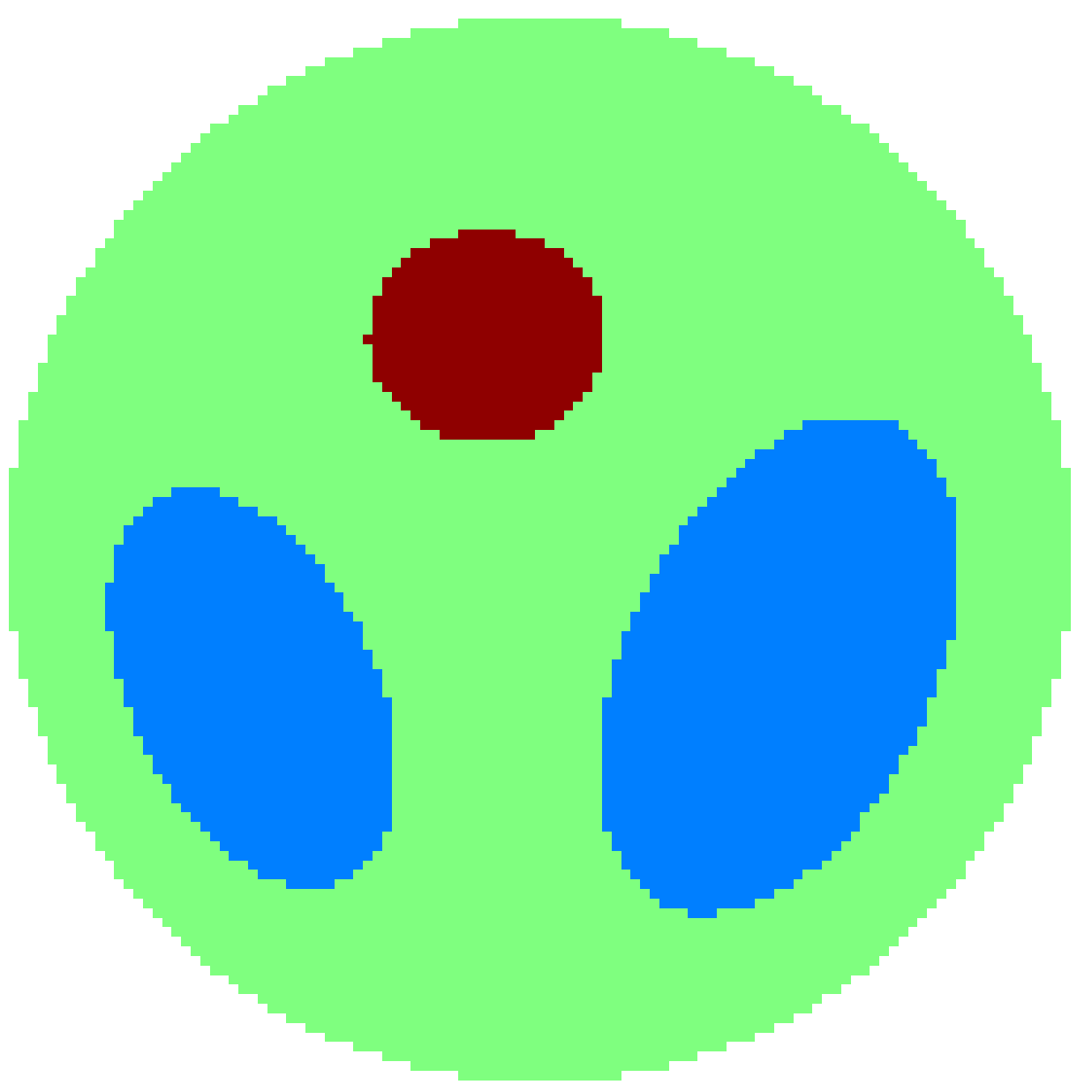}}
\put(70,2){\includegraphics[height=2.5cm]{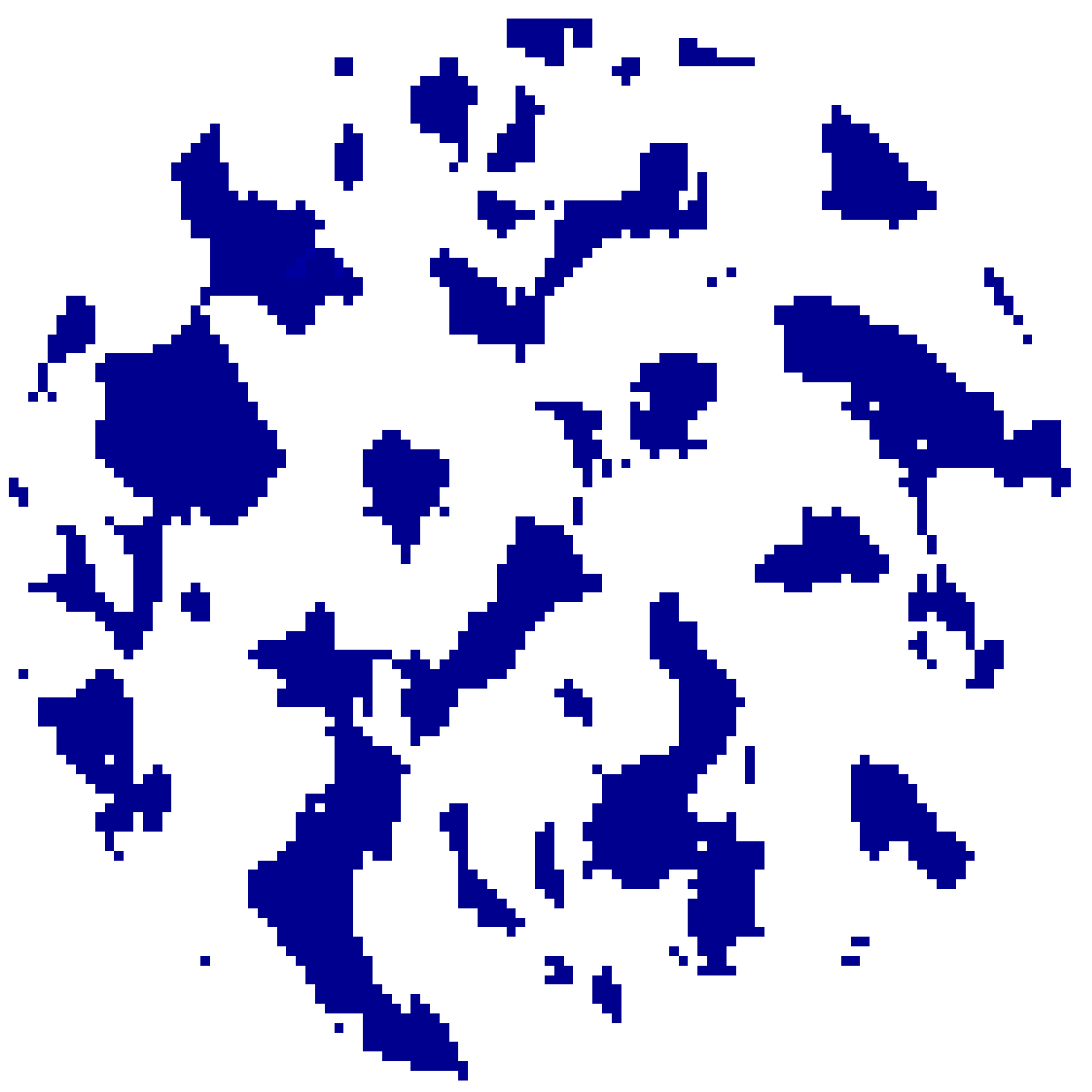}}
\put(170,2){\includegraphics[height=2.5cm]{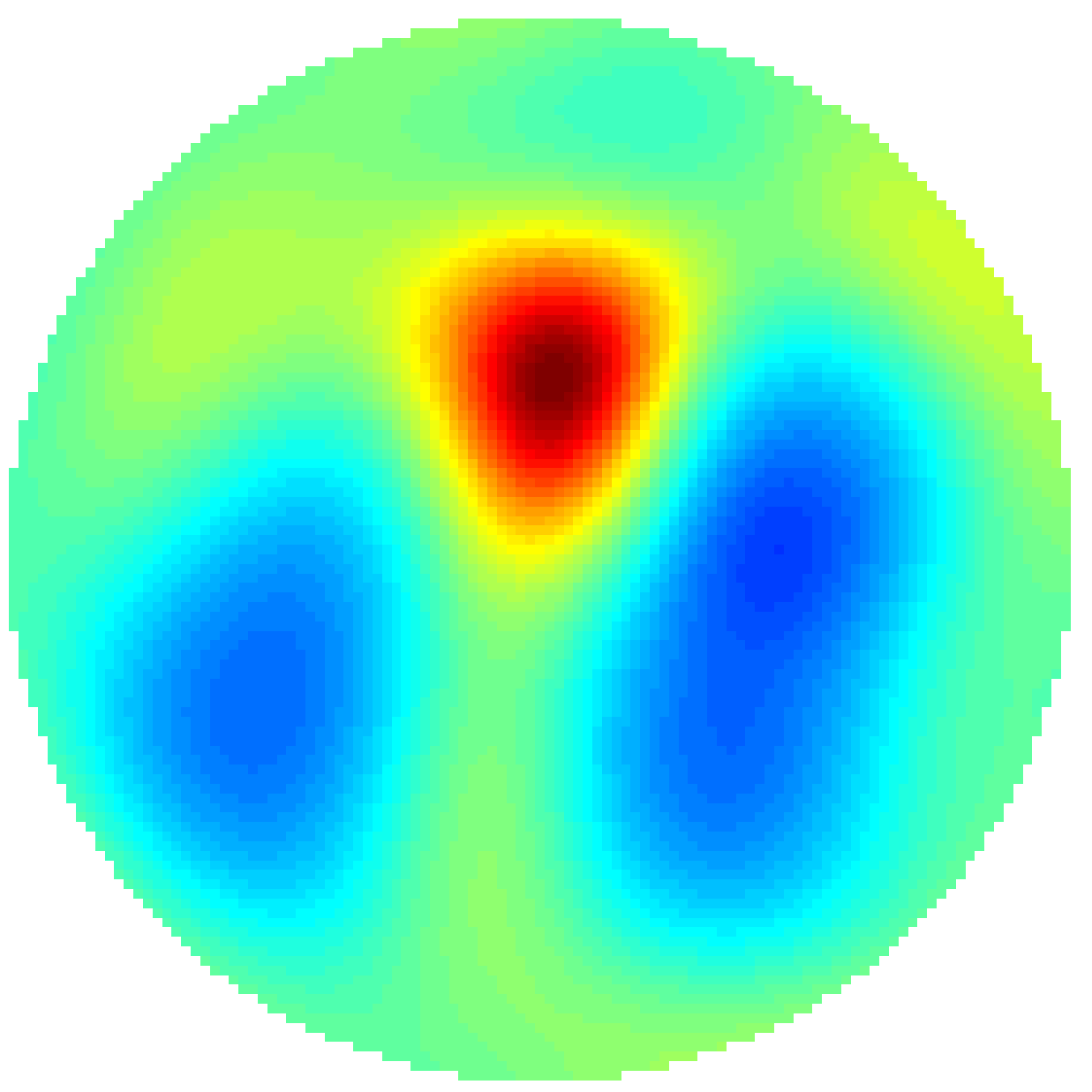}}
\put(270,2){\includegraphics[height=2.5cm]{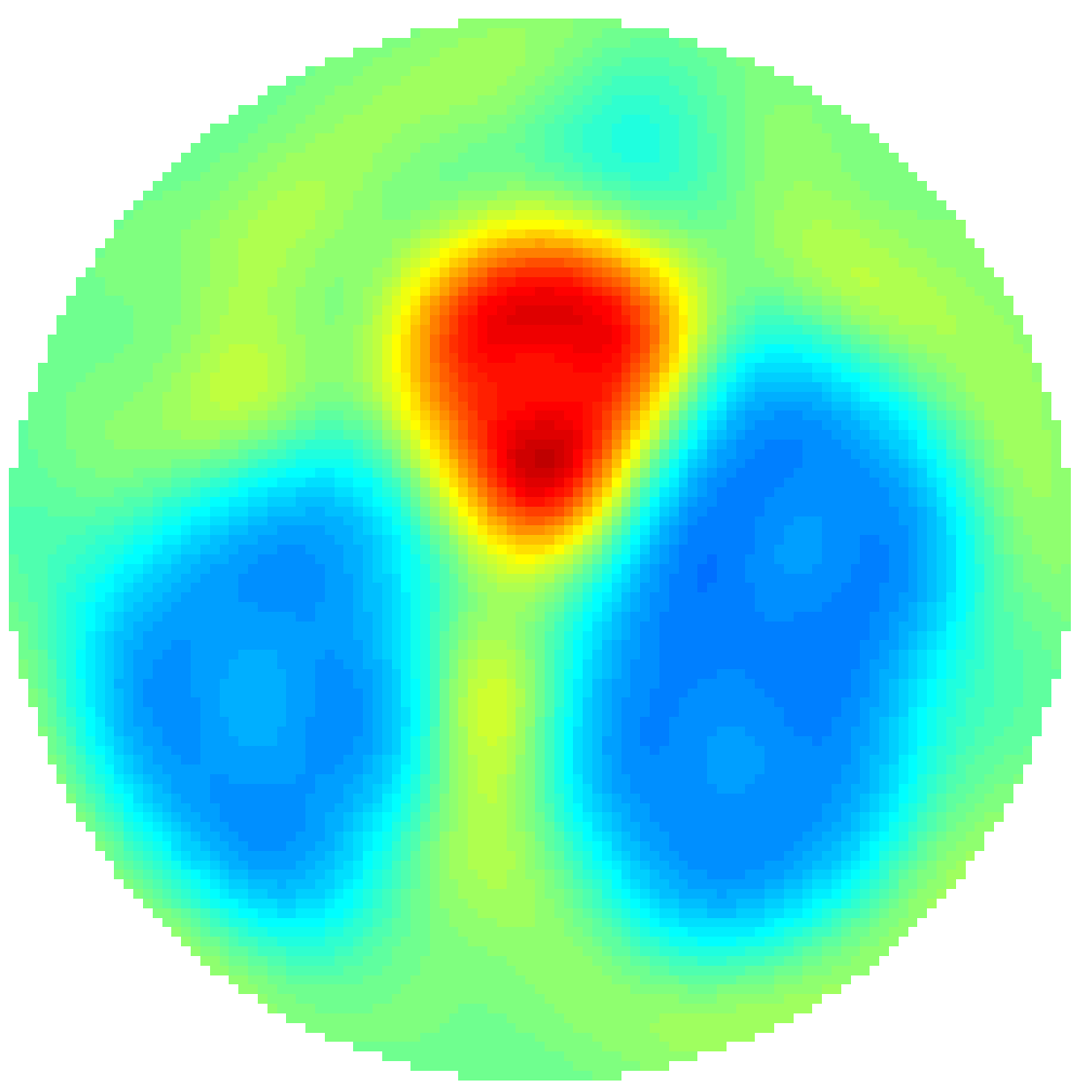}}
\put(-10,110){\includegraphics[height=1cm]{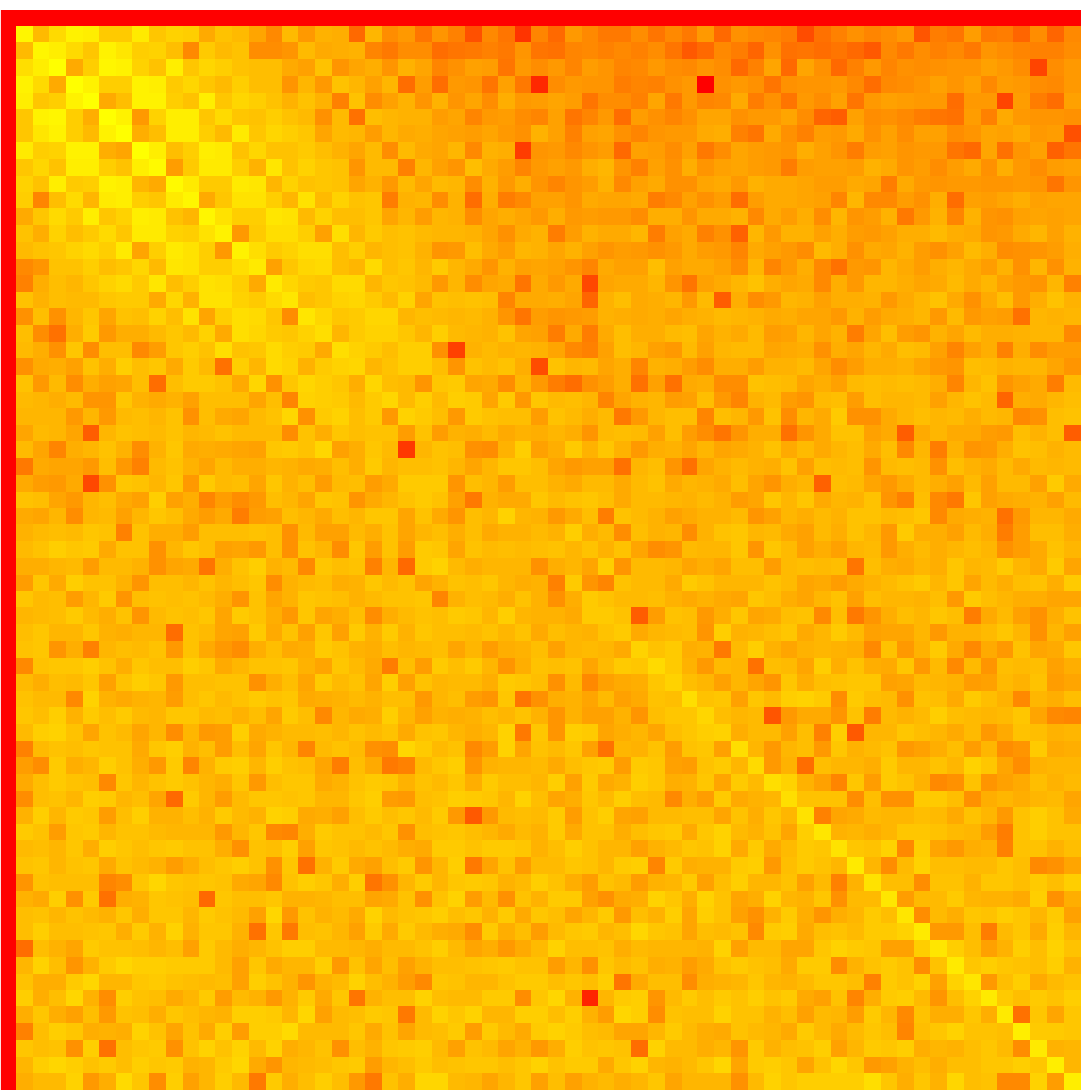}}
\put(70,160){\includegraphics[height=2.4cm]{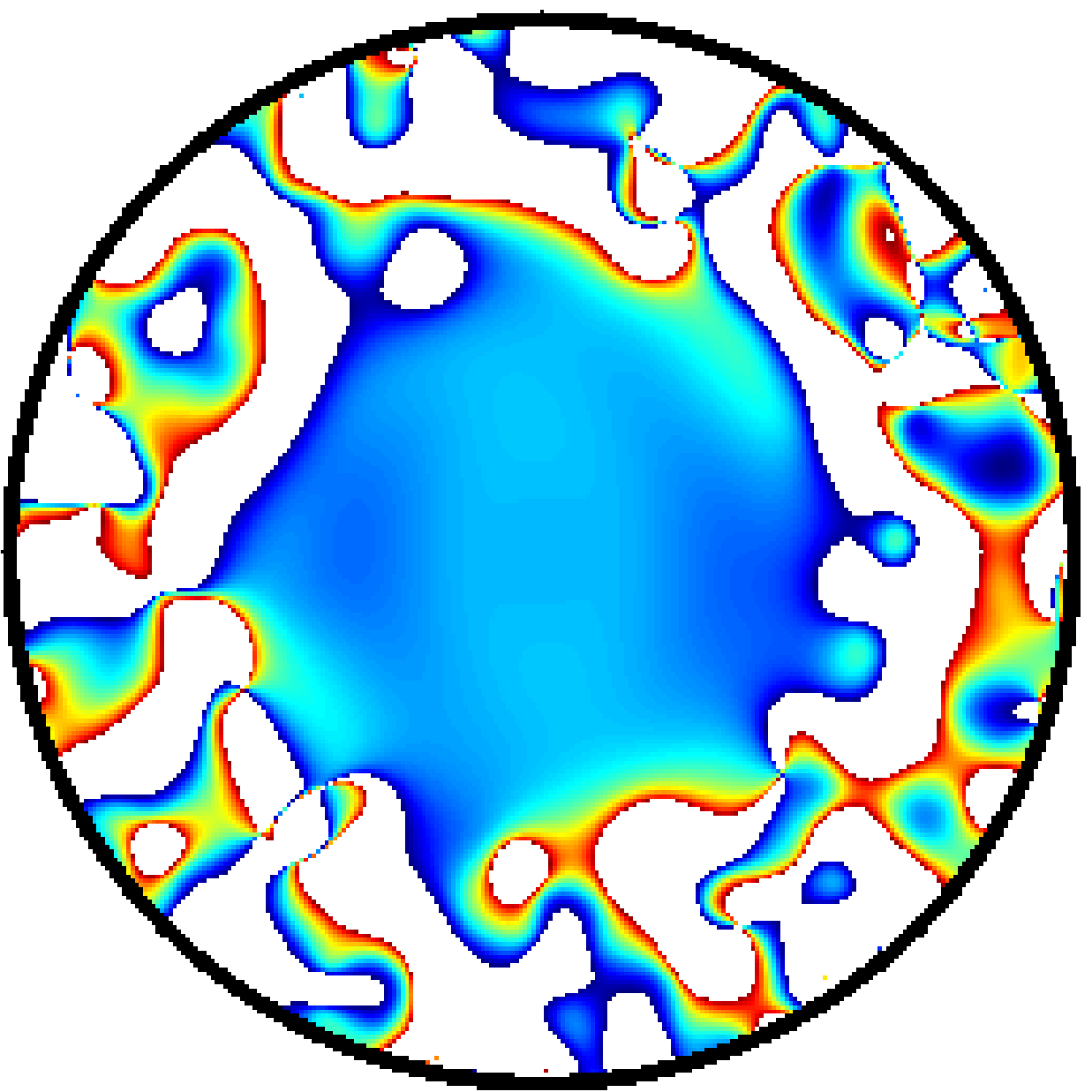}}
\put(170,160){\includegraphics[height=2.4cm]{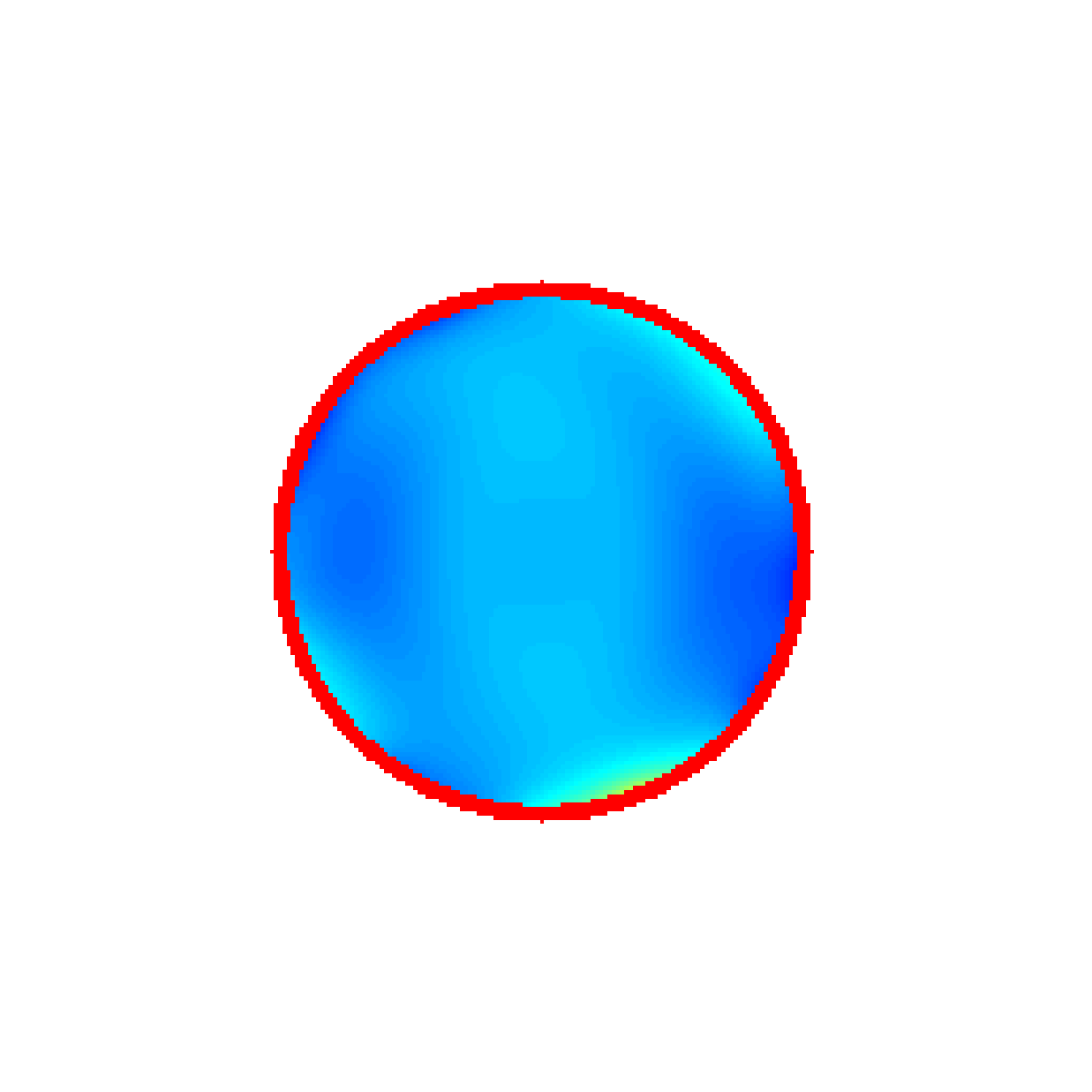}}
\put(270,160){\includegraphics[height=2.4cm]{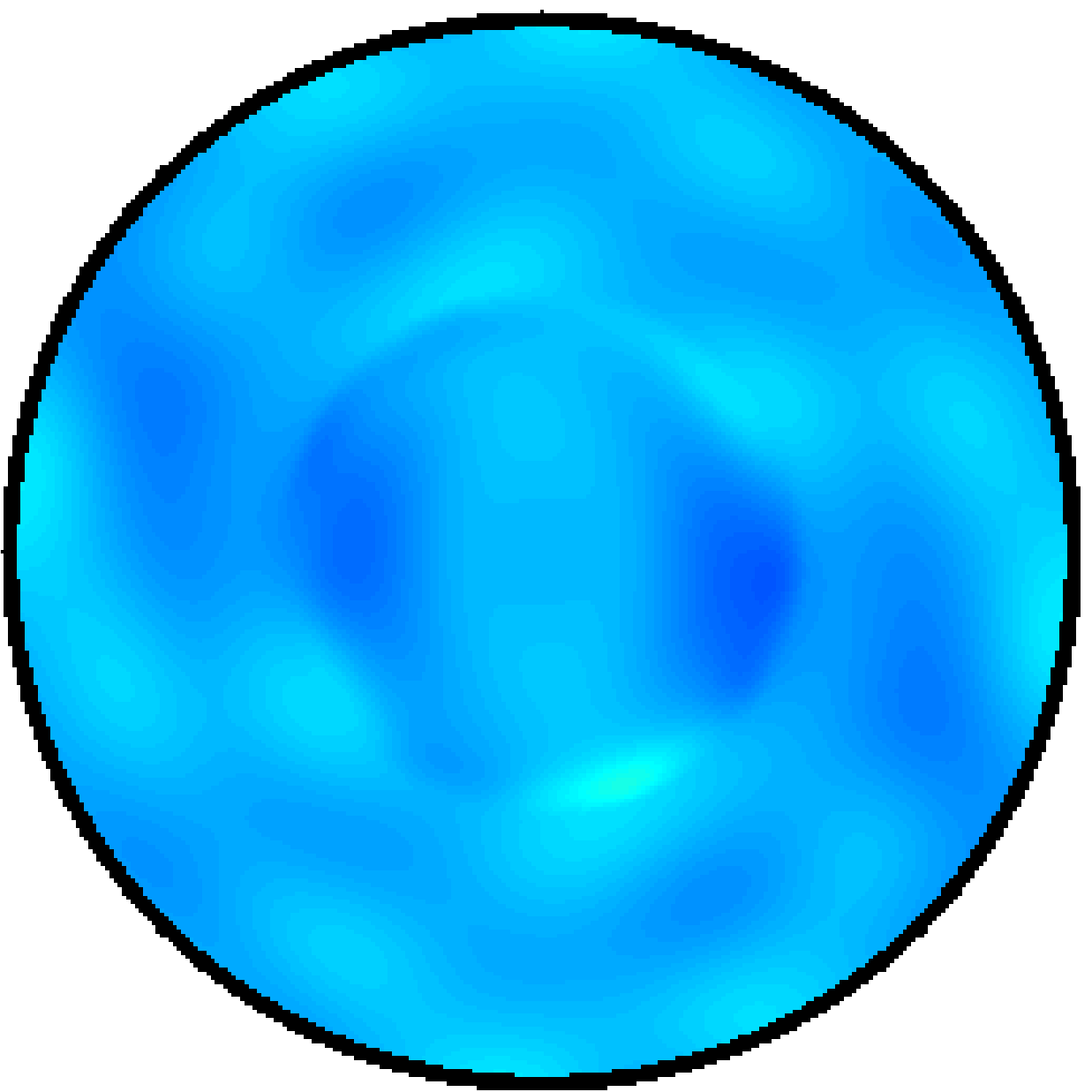}}
\thicklines \linethickness{.1mm} \put(5,75){\vector(0,1){28}}
\put(105,155){\vector(0,-1){80}} \put(205,155){\vector(0,-1){80}}
\put(305,155){\vector(0,-1){80}}
\put(309,140){\rotatebox{-90}{\footnotesize Inverse}}
\put(296,140){\rotatebox{-90}{\footnotesize transform}}
\put(209,140){\rotatebox{-90}{\footnotesize Inverse}}
\put(196,140){\rotatebox{-90}{\footnotesize transform}}
\put(109,140){\rotatebox{-90}{\footnotesize Inverse}}
\put(96,140){\rotatebox{-90}{\footnotesize transform}}
\put(22,140){\vector(1,1){40}}
\put(-27,120){$\Lambda_\sigma^\delta$}
\thinlines \multiput(-30,147)(4,0){97}{\line(1,0){2}}
\multiput(-30,90)(4,0){97}{\line(1,0){2}}
\put(215,79){\footnotesize Image domain}
\put(215,152){\footnotesize Frequency domain} \put(67,230){(a)}
\put(167,230){(b)} \put(267,230){(c)} \put(-33,0){(d)}
\put(67,0){(e)} \put(167,0){(f)} \put(267,0){(g)}
\end{picture}
\end{center}
\caption{\label{fig:scheme}Schematic illustration of the nonlinear
low-pass filtering approach to regularized EIT. The simulated
heart-and-lungs phantom (d) gives rise to a finite
voltage-to-current matrix $\Lambda_{\sigma}^\delta$ (orange
square), which can be used to approximately determine the
nonlinear Fourier transform (a). Measurement noise causes
numerical instabilities in the transform (irregular white patches
in (a)), leading to an unstable and inaccurate reconstruction (e).
However, multiplying the transform by the characteristic function
of the disc $|k|<5$ yields a lowpass-filtered transform (b), which
in turn gives a noise-robust approximate reconstruction (f). The
hybrid method presented in this paper uses {\em a priori}
information about the conductivity to estimate the missing part of
the nonlinear Fourier transform, resulting in (c). An improved
reconstruction (g) is achieved.}
\end{figure}

Generally speaking, regularization involves complementing insufficient measurement information by {\em a priori} knowledge about the conductivity. The D-bar method does this very explicitly assuming that the conductivity is twice continuously differentiable, which allows replacing the values of the nonlinear Fourier transform by zero in the unstable part of the frequency domain. Therefore, this low-pass filtering has the side effect that the resulting D-bar reconstructions are always smooth, as seen in the middle image of Figure \ref{fig:hammer}.

In many applications of EIT, including medical imaging, it is important to see boundaries between regions of different conductivities. The standard D-bar method practice of inserting zeroes for high frequencies is not ideal as crisp boundaries between different conductivity regions necessarily contain high frequencies.

\begin{figure}[t!]
\centering
\begin{picture}(330,470)
\definecolor{white}{gray}{1}
\thinlines

\put(0,400){\framebox(270,60)[lc]{}}
\put(10,442){\begin{minipage}[t]{8.5cm}
Set \textcolor{red}{$j=1$} and compute the {\bf nonlinear Fourier transform} $\T^{0}(k)$ in the stable disc $|k|<R$ from the noisy EIT data matrix $\Lambda_{\sigma}^{\delta}$.
\end{minipage}}

\put(50, 400){\vector(0,-1){40}}
\put(55,375){$\T^{0}$}

\put(0,300){\framebox(270,60)[lc]{}}
\put(10,342){\begin{minipage}[t]{8.5cm}
Compute the {\bf D-bar reconstruction} $\sigma^{(j)}_{\mbox{\tiny DB}}$ from the truncated transform $\T^{j-1}(k)$. If $j=1$ use cutoff disc $|k|<R$, else use larger disc $|k|<\widetilde{R}$.
\end{minipage}}

\put(50, 300){\vector(0,-1){40}}
\put(55,277){$\sigma^{(j)}_{\mbox{\tiny DB}}$}

\put(0,200){\framebox(270,60)[lc]{}}
\put(10,242){\begin{minipage}[t]{8.5cm}
Use TV segmentation flow to {\bf introduce edges} to the smooth image $\sigma^{(j)}_{\mbox{\tiny DB}}$. The result $\sigma^{(j)}_{\mbox{\tiny TV}}$ is a piecewise constant image.
\end{minipage}}

\put(50, 200){\vector(0,-1){40}}
\put(55,177){$\sigma^{(j)}_{\mbox{\tiny TV}}$}

\put(0,90){\framebox(270,70)[lc]{}}
\put(10,142){\begin{minipage}[t]{8.5cm}
{\bf Contrast enhancement.} Use the stable part of EIT data to modify $\sigma^{(j)}_{\mbox{\tiny TV}}$ into an optimal image called  $\sigma^{(j)}_{\mbox{\tiny CE}}$.
\textcolor{red}{If $j=J$ then return $\sigma^{(j)}_{\mbox{\tiny CE}}$ and stop.}
\end{minipage}}

\put(50, 90){\vector(0,-1){40}}
\put(55,67){$\sigma^{(j)}_{\mbox{\tiny CE}}$}

\put(0,-10){\framebox(270,60)[lc]{}}
\put(10,32){\begin{minipage}[t]{8.5cm}
Define $\T^{j}(k)$ for $R<|k|\leq \widetilde{R}$ as the new transform $\tilde{\T}^j(k)$ of $\sigma^{(j)}_{\mbox{\tiny CE}}$, for $|k|<R-1$ as $\T^{0}(k)$, and for $R-1<|k|\leq R$ as a blend of $\T^{0}$ and $\tilde{\T}^{j}$.
\textcolor{red}{Set $j:=j+1$.}
\end{minipage}}

\put(270, 30){\line(1,0){50}}
\put(320, 30){\line(0,1){275}}
\put(320, 305){\vector(-1,0){50}}
\put(283, 33){$\T^{j-1}$}
\put(283, 310){$\T^{j-1}$}

\put(-20, 450){(a)}
\put(-20, 350){(b)}
\put(-20, 250){(c)}
\put(-20, 150){(d)}
\put(-20, 50){(e)}

\end{picture}
\caption{\label{fig:flowchart}Flowchart of the proposed edge-preserving EIT reconstruction method.}
\end{figure}
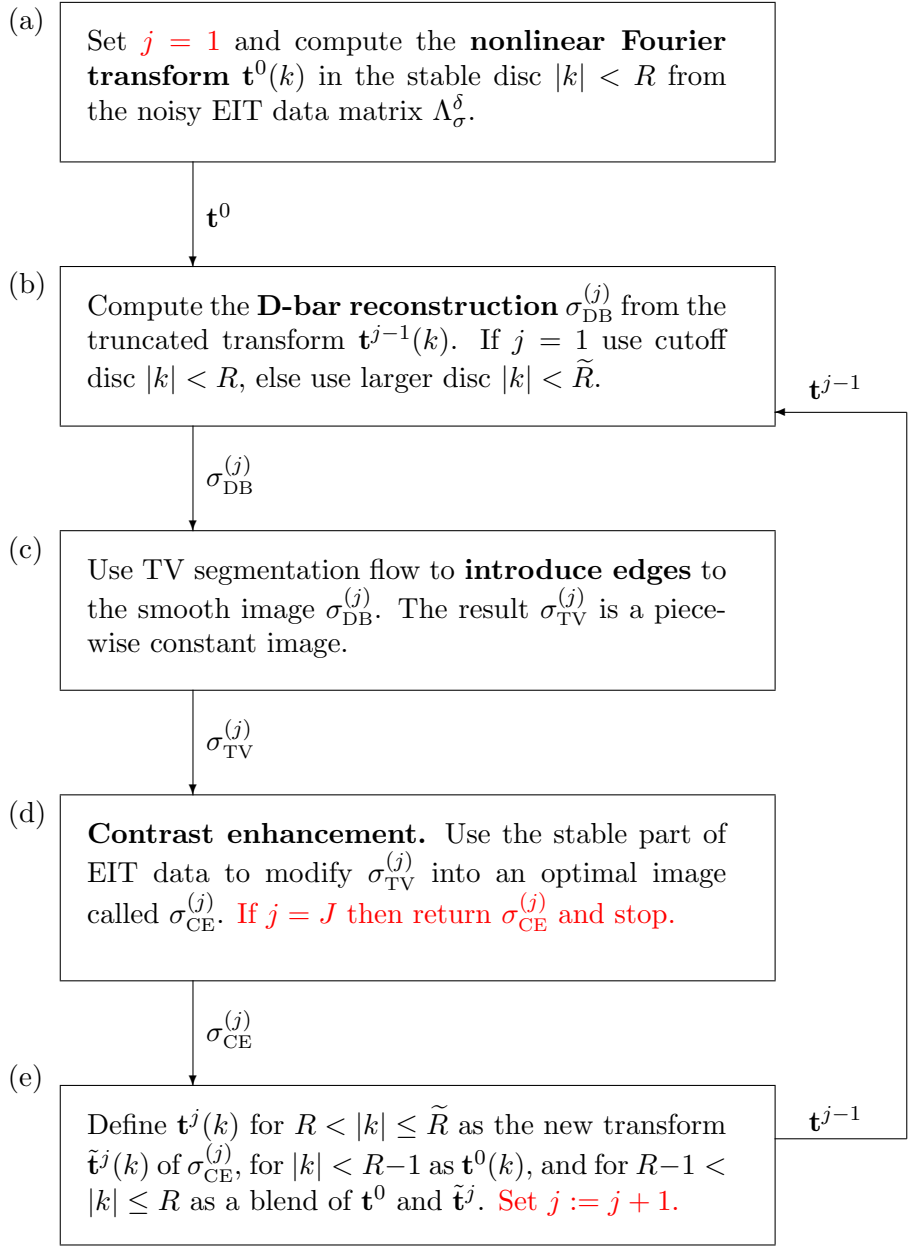

We introduce a novel edge-enhancing method for EIT, built upon the assumption that we know {\em a priori} that the conductivity is piecewise constant. It builds upon the stable D-bar reconstruction and increases the radius of reliable scattering data to pick up the missing high-frequency features (sharp edges and jumps) in the recovered conductivity.  The method applies Total Variation (TV) segmentation and data-driven contrast enhancement to the D-bar reconstruction regularized by low-pass filtering with cutoff frequency $R$. We exploit the methodology in \cite{Astala2006a,Astala2010,Astala2011} that allows the computation of nonlinear Fourier transform of the discontinuous segmented image in the annulus $R-1<|k|<\tilde{R}$, for certain $\tilde{R}>R$. This new transform is added on the annulus $R<|k|<\tilde{R}$ to the original transform restricted to the disc $|k|<R-1$, and we blend continuously both transforms on the annulus $R-1<|k|<R$. From this combined scattering data on the disc $|k|<\widetilde{R}$, a new sharper D-bar conductivity reconstruction is obtained. This procedure is iterated as outlined in Figure \ref{fig:flowchart}.

Figure \ref{fig:hammer} shows a reconstruction from simulated EIT data with $0.75\%$ relative noise added to the voltage data.  Our nonlinear method delivers a piecewise constant and edge-preserving reconstruction. The edges are more correctly located near the boundary. This is in accordance with the basic intuition about EIT: the deeper in you try to see, the harder it gets.

Let us comment on the variety of D-bar method we use. The three options are Schr\"odinger-type,
$2{\times}2$ (first-order) system type
and Beltrami-type.
They all have two steps: recover frequency-domain information and  reconstruct via an  inverse transform. Theoretically, only the Beltrami approach can deal with discontinuous conductivities. However, its second step is nonuniform in quality, see \cite[Section 6.3]{Astala2011} and \cite[Figure 12]{Astala2014}. We use below the {\em shortcut method} introduced in \cite{Astala2014}, combining Beltrami-type first step and Schr\"odinger-type second step.

The rest of this paper is organized as follows. In Section \ref{sec:model:review} we formulate the  mathematical EIT model and review the relevant literature. Section \ref{sec:Dbar} is devoted to a discussion of  the shortcut method.  The edge-promoting TV flow is described in Section \ref{sec:TVflow}.  Section \ref{sec:CE} introduces a contrast enhancement step based on the {\em CGO sinogram}, which is useful as a robust data-fidelity term. The new algorithm we present in this work is tested on two discontinuous phantoms for varying levels of noise. Section \ref{sec:numresults} outlines the computational details and in Section \ref{sec:numresults2} the numerical results are presented and discussed. We conclude our findings in Section \ref{sec:concl}.

\section{Mathematical model and literature review}\label{sec:model:review}

\noindent
We concentrate on the two-dimensional case with $\Omega$ representing the unit disc. However, all our techniques can be extended to simply connected domains $\Omega\subset\R^2$ with Lipschitz boundary $\bndry$.   Throughout the paper  the following notation is used.  For $r>0$, $D(0,r)$ denotes the disc in
the plane centered at the origin with radius $r$. For any open set
$U$ in the plane, $\overline{U}$ denotes the closure of $U$. We associate $\C$ and $\R^2$ by $z=(x,y)=x+iy$.


Let $f$ denote the electric voltage potential maintained at the boundary.
The corresponding potential $u$ inside the domain $\Omega$ satisfies the Dirichlet problem for the elliptic conductivity equation
\begin{equation}\label{eq-conduct}
\begin{array}{rcl}
\nabla\cdot\sigma(z)\nabla u(z) & = &0,\qquad z\in\Omega\\
u|_{z\in\bndry} &=& f(z),\quad z\in\bndry,
\end{array}
\end{equation}
where $\sigma\in L^\infty(\Omega)$ is an isotropic conductivity satisfying $\sigma(z)\geq c>0$.  The \emph{Dirichlet-to-Neumann} (D-N) map, is defined by
\begin{equation}\label{eq-DN-def}
\Lambda_\sigma f=\left.\sigma\frac{\p u}{\p\nu}\right|_{\bndry},
\end{equation}
where $\nu$ denotes the outward facing unit normal vector to the boundary. It is well-known that $\Lambda_\sigma:H^{1/2}(\bndry)\rightarrow H^{-1/2}(\bndry)$ is a bounded linear operator. One can think of $\Lambda_\sigma$ as a mathematical model for voltage-to-current measurements performed at the boundary.

The above mathematical model for the inverse conductivity problem was first formulated by Alberto  Calder\'{o}n in 1980 \cite{Calder'on1980}. Inspired by his engineering background, he asked whether it is possible to calculate a bounded conductivity $\sigma(z)\geq s>0\in L^\infty(\Omega)$ in terms of electrical boundary measurements. He gave a proof for the linearized problem in the case of infinite-precision knowledge of the noise-free D-N map $\Lambda_\sigma$.

In practical EIT, one needs to recover the electric conductivity distribution $\sigma:\Omega\ra\R$ in a regularized manner  from noisy boundary measurements $\Lambda_\sigma^\delta$, where $\|\Lambda_\sigma-\Lambda_\sigma^\delta\|_Y\leq \delta$ for a known noise level $\delta$ and an appropriate norm $\norm{\cdot}_Y$. Among currently available EIT algorithms, the D-bar method is the only one with a regularization analysis  \cite{Knudsen2009}.

The D-bar method is based on a nonlinear Fourier transform, which is not physically measurable. The definition of the transform depends on certain  ``almost exponential'' functions, the so-called {\em complex geometric optics} (CGO) solutions. CGO solutions were introduced by Faddeev in 1966 \cite{Faddeev1966} and later introduced in the context of inverse problems by Sylvester and Uhlmann in 1987 \cite{Sylvester1987}.

The D-bar method was invented by Beals and Coifman for use in nonlinear evolution equations the early 1980's, e.g. \cite{Beals1985}. The first description of the use of the D-bar approach in inverse problems in dimensions $n\geq 2$ was published by R.G. Novikov in Russian in 1987 (for English translation see \cite{Novikov1988}). The Schr\"odinger equation approach is still used in dimensions three and higher. Rigorous mathematical theory of two-dimensional D-bar methods fall into the following three categories.

{\bf Schr\"odinger-type approach.} The first uniqueness result for the inverse conductivity problem in dimension two was published by Nachman in 1996 for $\sigma\in W^{2,p}$, $p>1$, in \cite{Nachman1996}. The first computational implementation of a D-bar method, based on a kind of Born approximation, was published in 2000 \cite{Siltanen2000}. The full nonlinear algorithm was developed and equipped with a regularization strategy in \cite{Mueller2003,Knudsen2004,Knudsen2007,Knudsen2009}. For reconstructions from experimental data, see \cite{Isaacson2006,Murphy2007a,DeAngelo2010,Dodd2014,Herrera2015}.	

{\bf The $2{\times}2$ system approach.}
The smoothness assumption for the conductivity was reduced by Brown and Uhlmann to $\sigma\in W^{1,p}$, $p>2$, see \cite{Brown1997}. For further theory and computational implementation, see \cite{Knudsen2004a,Knudsen2002}. This approach has the benefit of being applicable to complex-valued $\sigma$ (conductivity$+i\cdot$permittivity), as shown theoretically by Francini \cite{Francini2000}. For computational results, see \cite{Hamilton2012,Hamilton2013,Herrera2015}.

{\bf Beltrami-type approach.}
Astala and P\"{a}iv\"{a}rinta introduced a theoretical reconstruction method for $\sigma\in L^\infty$ conductivities in \cite{Astala2006, Astala2006a}. A corresponding computational reconstruction method was developed in \cite{Astala2010,Astala2011} and further analysed in \cite{Astala2014}.

Summarizing, numerical D-bar reconstruction methods provide direct, non-iterative, parallelizable, nonlinear  reconstruction approaches which are effective on experimental EIT data and even in real-time applications \cite{Dodd2014}. There are also three-dimensional D-bar type methods \cite{Novikov1988,Nachman1988,Cornean2005,Boverman2008a,Bikowski2011,Delbary2014}.

Conditional stability estimates (involving non-noisy data) for these D-bar methods appeared in \cite{Liu1997,Barcelo2001,Barcelo2007,Clop2010,Beretta2010,Faraco2013} and in dimension greater than two, the most recent stability result is in \cite{Caro2013}. For regularization analysis involving noisy data see \cite{Knudsen2009}.

Several edge-preserving regularization methods have been suggested for EIT in the literature, including Total Variation and sparsity-promo\-ting techniques \cite{Dobson1994,Kaipio2000,Rondi2001,Chung2005,vandenDoel2006,Jin2012a}. The method proposed in this paper is very different from all of these approaches as it uses a proven regularized D-bar image as its starting point and futher uses the inverse scattering of the D-bar methodology to guide the image segmentation.  The proposed method differs from \cite{Hamilton2014}, where the {\em CGO sinogram} was introduced, by instead focusing on enlarging the scattering radius to produce more accurate conductivity reconstructions.

\section{The ``shortcut'' D-bar reconstruction method}\label{sec:Dbar}

\subsection{The Schr\"odinger-type D-bar method}

\noindent
The regularized D-bar method \cite{Nachman1996,Siltanen2000,Knudsen2009} for $C^2$-conductivities  consists of two steps, namely $\Lambda^\delta_\sigma \overset{1}{\longrightarrow} \T_R(k)\overset{2}{\longrightarrow} \sigma_R(z)$. Step 1 is not used here; we refer the reader to \cite[Chapter 15]{Mueller2012} for details.

Step 2 goes from truncated scattering data $\T_R:\C\ra \C$, supported in the disc $|k|<R$, to the regularized conductivity as follows. For each $z\in\Omega$, solve the integral equation
\begin{equation}\label{eq-dbark-mu-sol}
m_{R}(z,k)=1+\frac{1}{(2\pi)^2}\int_{|\kappa|<R}
\frac{\T_R(\kappa)}{(k-\kappa)\bar{\kappa}}
e(-z,\kappa)\,\overline{m_{R}(z,\kappa)}\;d\kappa_1 d\kappa_2,
\end{equation}
where $e(z,k):=\exp\left\{i\left(kz+\bar{k}\bar{z}\right)\right\}=\exp\{2i\,\mbox{Re}(kz)\}$.  The regularized conductivity is computed by
\begin{equation}\label{eq-muR-to-sigR}
\sigma_R(z)=\left(m_{R}(z,0)\right)^2.
\end{equation}

\subsection{The Beltrami-type D-bar method}

\noindent The D-bar method  for $L^\infty$-conductivities \cite{Astala2006, Astala2006a,Astala2011} is based on complex geometric optics solutions $f_{\pm\mu}(z,k)$ to the Beltrami equation
\begin{equation}\label{eq:BeltramiEq}
\overline{\partial}_z f_{\pm\mu}(z,k)=\pm\mu(z,k)\overline{\dz f_{\pm\mu}(z,k)},
\end{equation}
where $\mu(z)=(1-\sigma(z))/(1+\sigma(z))$
and $f_{\pm\mu}(z,k)=e^{ikz}(1+\mathcal{O}(\frac{1}{|z|}))$ as $|z|\ra\infty$.  
Set $M_{\pm\mu}(z,k)=e^{-ikz}f_{\pm\mu}(z,k)$.
The reconstruction method has three steps:
\[\Lambda^\delta_\sigma \overset{1}{\longrightarrow} M_{\pm\mu}(\cdot,k)\vert_{\bndry}\overset{2}{\longrightarrow}\begin{array}{c}
\text{\small{Transport}}\\
\text{\small{Matrix}}
\end{array}\overset{3}{\longrightarrow} \sigma_R(z),\]
out of which we only use Step 1 (for Steps 2 and 3, see \cite[Section 16.3]{Mueller2012}). For each fixed $k\in\C$, $|k|<R$, solve
\begin{equation}\label{eq:CGOBeltrami_traces_Mmu}
M_{\pm\mu}(\cdot,k)\vert_{\bndry}+1=\left(\mathcal{P}^k_{\pm\mu}+\mathcal{P}_0\right)M_{\pm\mu}(\cdot,k)\vert_{\bndry}
\end{equation}
to obtain the CGO traces $M_{\pm\mu}(\cdot,k)$ for $z\in\bndry$, where $\mathcal{P}^k_{\pm\mu}$ and $\mathcal{P}_0$ are the projection operators described in \cite{Astala2011}.

\subsection{The shortcut method}\label{sec:shortcutmethod}

There is a connection \cite[Section 5]{Astala2010} between the Beltrami CGOs of
\cite{Astala2006,Astala2006a} and Schr\"{o}dinger CGOs of
\cite{Nachman1996}, as well as their associated scattering
transforms $\tau(k)$ and $\T(k)$, respectively \cite[Section 6]{Astala2010}. Namely, when  $\sigma\in\C^2$ we have
\begin{equation}\label{eq-tau-to-t}
  \T(k)
  =
  -4\pi i\bar{k}\tau(k)
  =
 -2 i\bar{k} \int_{\R^2}\overline{\left(\overline{\partial}_z \middle(M_{+\mu}(z,k) - M_{-\mu}(z,k) \right)}
dz_1 dz_2.
\end{equation}
Numerical evidence \cite{Astala2014} suggests that the above connection holds for $\sigma$ with jump discontinuities, resulting in numerical equivalence between the reconstructed conductivities.  In fact, the solution of \eqref{eq-dbark-mu-sol} is faster and more stable than the transport matrix method of \cite{Astala2011}.
Therefore, we will use the combined D-bar algorithm presented in \cite{Astala2014}, called the {\it shortcut method}, which has the following steps:
\begin{itemize}
\item[\textbf{Step 1:}] From noisy boundary measurements $\Lambda^\delta_\sigma$ to CGO boundary traces $M_{\pm\mu}(\cdot,k)\vert_{\bndry}$ by solving \eqref{eq:CGOBeltrami_traces_Mmu}.

\item[\textbf{Step 2:}] From boundary CGO traces $M_{\pm\mu}(\cdot,k)\vert_{\bndry}$ to truncated scattering data $\T_R(k)$.  The analyticity of $M_{\pm\mu}(\cdot , k)$ outside $\Omega$ leads to the following development for $|z|>1$
\[M_{\pm\mu}(z , k) = 1 + {a_1^\pm (k)\over z} + {a_2^\pm (k)\over z^2} + \ldots\]
from which one defines $\T_R(k) := -4\pi i \bar{k} \tau_R(k) $, where
\begin{equation}\label{eq:tau_Beltrami_from_Mu}
  \tau_R(k) := \begin{cases}
 {1\over 2}\left(\overline{a_1^+(k)}-\overline{a_1^-(k)}\right), & |k|<R\\
0 & |k|\geq R.
\end{cases}
\end{equation}

\item[\textbf{Step 3:}] From truncated scattering data $\T_R(k)$ to conductivity $\sigma_R(z)$. Solve \eqref{eq-dbark-mu-sol} and evaluate \eqref{eq-muR-to-sigR}.
\end{itemize}



\section{Image segmentation method}\label{sec:TVflow}

\noindent
In this section, we present a Mumford-Shah based segmentation model.  For the simplicity of notation, let $\sigma_0(x)$ be the input conductivity  image defined on $\Omega$.  The image segmentation problem is to find a partition of $\Omega$ into $K$ disjoint subdomains $\{\Omega_k\}_{k=1}^K$, i.e.
\[
\Omega=\bigcup_{k=1}^{K}\Omega_{k}\, ; \qquad \Omega_{k}\cap\Omega_{j}=\emptyset  \textrm{ for }\, k\neq j.
\]

The celebrated Mumford-Shah \cite{mumford1989optimal} variational segmentation problem is as follows:  find a piecewise smooth function and  a partition edge (closed) set $\Gamma=\bigcup_{k=1}^K \partial\Omega_k$ such that the following functional is minimized:
\[
E(f,\Gamma) = \frac{\lambda}{2}\int_\Omega(f(x)-\sigma_0(x))^2 dx + \int_{\Omega\backslash \Gamma}|\nabla f(x)|^2{dx}  +\alpha \mathcal{H}^{1}(\Gamma),
\]
where  $\mathcal{H}^1(\Gamma)$ denotes the 1D Hausdorff measure of the edge set $\Gamma$.  Due to the  complexity of the model, many simplifications are proposed. In the simplest form of the model, $f$ is assumed to be piecewise constant on each $\Omega_k$ and the model is reduced to
\begin{equation}
\label{eq:msmodel}
\min\limits_{\{\Omega_k, c_k\}_{k=1}^K} \left\{\frac{\lambda}{2} \sum_{k=1}^K \int_{\Omega_k} (\sigma_0(x)-c_k)^2dx+\sum_{k=1}^K |\partial\Omega_k| \right\},
\end{equation}
where $c_k\in\bfR$ for $k=1,\dots, K$ is the mean intensity  for each subregion $\Omega_k$. The parameter $\lambda>0$ is used to balance the data fitting and the total length of regions interfaces. This model \eqref{eq:msmodel} is  hard to be solved directly.  In fact, when the regions $\Omega_k$ are determined, the optimal $c_k$ is given as
 \begin{equation}
 \label{eq:eqc}
  c_k=\frac{\int_{x\in\Omega_k} \sigma_0(x)}{|\Omega_k|}.
 \end{equation}

Thus  we can consider an alternating scheme on solving $\Omega_k$ and $c_k$ iteratively. Once $c_k$ is determined, we solve for the $\Omega_k$. By doing so, we  introduce the labeling function $u_k$ of the  disjoint subregions $\Omega_k$
\begin{equation}
u_k(x) = \left\{ \begin{array}{ll}
1 & \textrm{if $x\in\Omega_k$}\\
0 & \textrm{otherwise}\\
\end{array} \right., \quad\mbox{ for } k=1, \cdots, K.
\end{equation}
According to the co-area formula, the perimeter of a set $\Omega_k$ is given by the total variation of $u_k$
\begin{equation}
\label{eq:coarea}
|\partial\Omega_k|=\int_\Omega |D u_k|
\end{equation}
and the total variation $\int_\Omega |D u|$ is defined in distribution sense
\begin{equation}
\label{eq:TVdual}
\int_\Omega |Du|:=\sup\Big\{-\!\int_\Omega \! u \,\mbox{div}\phi\, dx: \phi\in C_c^\infty(\Omega;\bfR^d), |\phi(x)|\leq 1,\mbox{ a.e. } x\in\Omega\Big\}.
\end{equation}
It is well-known that if $u\in W^{1,1}(\Omega)$ then  $\int_\Omega |D u|=\int_\Omega |\nabla u(x)|dx$.

Meanwhile, since each pixel can be only assigned to be one region, the labeling function $u_k(x)$ satisfies the following constraint:
\begin{equation}
  \sum_{k=1}^K u_k(x)=1,\mbox{ a.e. } x\in\Omega
\end{equation}

 Generally, convex relaxation is made by allowing $u_i$ to take values continuously in $[0,1]$ to overcome the computation complexity of the binary constraint.
The overall  model is reformulated as
\begin{equation}
\label{eq:genericnonconvex}
\begin{aligned}
&\min_{\scriptsize{\{u_k, c_k\}_{k=1}^K}} \left\{\sum_{k=1}^K \int_\Omega|D u_k|+\sum_{k=1}^K \int_{\Omega} u_k f_k\right\}\\
s. t.&\quad (u_1(x), \cdots, u_K(x))\in \bfS
\end{aligned}
\end{equation}
where we denote
\begin{equation}
\label{eq:fk}
f_k(x)=\frac{\lambda}{2} (\sigma_0(x)-c_k)^2,
\end{equation}
the constraint set
\begin{equation}
\label{eq:simplex}
\bfS=\left\{ (u_1,\cdots, u_K)\in \mbox{BV}(\Omega, [0,1]^K): \  \sum_{k=1}^K u_k(x)=1,\mbox{ a.e. } x\in\Omega\right\}
\end{equation}
and  $\mbox{BV}(\Omega,[0, 1]^K)$ denotes the bounded variation functions  product space valued on $[0,1]^K$.

With fixed $c_k$, the  convex relaxed formulation for $u_k$ allows us to develop efficient algorithms based on the well studied total variation minimization.  For example, it has been extensively studied in \cite{Chan2005, Chan06,Bae2011,Pock09,Chambolle2011}. Theoretically,  the global solution to the original binary model can be achieved for the case of two regions when the intensities $c_k$ are given.

This above region-based model can be further  combined with an edge based approach to improve the segmentation quality and speed,  such as in \cite{sandberg2005logic,Bresson2007,Bae2011}.  Assuming $u_k\in W^{1,1}(\Omega)$, the weighted total variation model can be defined  as
\begin{equation}
\label{eq:convex_edge_tv} J(\bfu)=\sum_{k=1}^K
\int_{\Omega}g(x)|\nabla u_k(x)|dx
 \end{equation}
where $g(x)\geq 0$ is an edge function taking small values at locations with large gradient and large values for smooth region.  For example,  a usual choice is
\begin{equation}
\label{eg:gx}
g(x)=\frac{1}{1+s \|\nabla \tilde{\sigma_0}(x)\|^2}
\end{equation}
where $\tilde{\sigma_0}$ is a smoothed version of the given image $\sigma_0$ and $s>0$ is a positive number. Note that if $g(x)$ is identical to $1$, it reduces to the model \eqref{eq:genericconvex}.

In the following, we present a primal-dual splitting method used in \cite{Pock09,EZC10,chambolle2011first,TZS2014}. Denote $\bfu=(u_1, \cdots,u_K)$, $\bff=(f_1,\cdots, f_K)$ and $J(\bfu)= \sum_{k=1}^K
\int_{\Omega}g(x)|\nabla u_k|dx$ and $\langle\bfu, \bff\rangle=\sum_{k=1}^K \int_{\Omega} u_k(x)f_k(x)dx$, then
the  convex minimization problem is rewritten as:
\begin{equation}
\label{eq:genericconvex}
\bfu^*=\arg\min_{\scriptsize\bfu\in\bfS}\left\{ J(\bfu)+\langle\bfu, \bff\rangle\right\}
\end{equation}

Based on the dual definition \eqref{eq:TVdual}, we consider the following min-max model:
\begin{equation*}
\min_{{\scriptsize{\bfu\in\bfS}}}\max_{{\scriptsize{\bfp\in\bfT}}} E(\bfu,\bfp)=\left\{ \langle\bfu, \mbox{\scriptsize{div}}(\bfp)\rangle+\langle\bfu, \bff\rangle\right\}
\end{equation*}
where $\bfp=(p_1,\cdots, p_K)$ for $p_k(x)\in p_k\in C_c^\infty(\Omega;\bfR^d)$, $\langle\bfu, \mbox{\scriptsize{div}}(\bfp)\rangle=\sum_{k=1}^K \int_\Omega u_k(x)\mbox{div}p_k(x)dx$
 and
\begin{eqnarray*}
  \bfT
  &=&
  \Big\{\bfp=(p_1,\cdots, p_K)\,:\,  \Big(\sum_{i=1}^d|p_k^d(x)|^2\Big)^{1/2}\leq g(x),\mbox{ a.e. } x\in\Omega \\
  &&
  \mbox{ and for all } k=1,\cdots, K \Big\}.
\end{eqnarray*}

 The specific algorithm is given as follows:
\begin{framed}
\begin{itemize}
\item[\textbf{Step 0:}]  \textbf{Initialization:} Choose $c_1,\cdots, c_K$ as the initial guess of the mean intensity of each region.  $\tau_1,\tau_2> 0$ are parameters such that $\tau_1\tau_2\leq 1/8$ and $(\bfp^0, \bfu^0)\in \bfT\times \bfS$. $\overline{\bfu}^0=\bfu^0$.

\textit{Set $i:=0$ and run the outer loop  as follows:}

\vspace{1em}
\item[\textbf{Step 1:}]
 \textit{Set $j:=0$ and run the inner loop  to compute $u_k(x)$ for $k=1,\cdots, K$ and $x\in\Omega$}
 \begin{itemize}
\item[\textbf{Step 1.1:}]

Compute the dual variable $\bfp^{j+1}=\Pi_{\bfT}(\bfp^j+\tau_1 \nabla\overline{\bfu}^j)$
where $\Pi_{\bfT}(\cdot)$ denotes the projection operator onto the convex set $\bfT$.

\vspace{0.5em}

\item[\textbf{Step 1.2:}]   Compute the primal variable $\bfu^{j+1}=\Pi_{\bfS}(\bfu^j+\tau_2(\bff_i+\mbox{div}^*\bfp^{j+1}))$
where $\Pi_{\bfS}(\cdot)$ denotes the projection operator onto the convex set $\bfS$.

\vspace{0.5em}

\item[\textbf{Step 1.3:}]  Compute the auxiliary primal variable: $\overline{\bfu}^{j+1}=\bfu^{j+1}+(\bfu^{j+1}-\bfu^j)$

\vspace{0.5em}

\item[\textbf{ Step 1.4:}] Set $j=j+1$ update until stopping conditions satisfied, output $\bfu$.

\vspace{1em}
\end{itemize}
\item [\textbf{Step 2:}] Compute the piecewise regions $\Omega_k$ by the binarization of $u_k(x)$.  Generally, a global minimizer of \eqref{eq:genericconvex} might not be binary, and a final thresholding step needs to be taken to get a binary solution
\begin{equation}
\label{eq: finalshrink}
    u_k^*(x) = \left\{
    \begin{array}{ll}
    1 & \textrm{if $u_k^*(x)= \max\{u_1(x), u_2(x), \ldots, u_K(x)\}$} \\
    0 & \textrm{otherwise} \\
    \end{array}  \right.
\end{equation}
If the maximizer is not unique, the maximizer with smallest subscript is  used as a convention.
\item [\textbf{Step 3:}] Update the mean intensity estimation: $c_k^i$ for $k=1, \cdots, K$ by \eqref{eq:eqc} and $\bff_i$ by \eqref{eq:fk}.

\vspace{1em}

\item [\textbf{Step 4:}] Set $i=i+1$ update until stopping conditions satisfied.

\vspace{1em}

\end{itemize}
\end{framed}
  After we obtain the label functions $u^*_k(x)$ and $c_k$ for $k=1, \cdots, K$, an image can be reconstructed as a piecewise constant function with the mean intensity $c_k$ in the corresponding $k-$ th region, i.e.
 \begin{equation}
 \label{eq:segregularized}
 \sigma_{TV}(x)= \sum_{k=1}^K u_k^*(x)c_k
 \end{equation}
 Hence, this reconstructed image can be used as  a piecewise constant regularized approximation to the original image by letting $\sigma_0(x):=\sigma_{DB}(x)$ in the whole algorithm described in Figure \ref{fig:flowchart}.


\section{Data-Driven Contrast Adjustments}\label{sec:CE}

\subsection{The Beltrami CGO sinogram}\label{sec:CGOsino}

\noindent
We extend the concept of the CGO sinogram, introduced in \cite{Hamilton2014}, to discontinuous conductivities. Set
$$
  \mathcal{S}_\sigma(\theta,\varphi,\rho):=
M_{\mu}(e^{i\theta},\rho \, e^{i\varphi})-1,
$$
where $z=e^{i\theta}$ and $k=\rho e^{i\varphi}$ for $\theta,\varphi\in [-\pi,\pi)$ and the traces $M_{\mu}(e^{i\theta},\rho \, e^{i\varphi})$ of the CGO solutions are solved from the noisy EIT data $\Lambda_{\sigma}^{\delta}$ using equation (\ref{eq:CGOBeltrami_traces_Mmu}). The radius $\rho$ must be smaller than the noise-dependent cutoff frequency $R$.

The traditional data-fidelity term used in EIT is $\|\Lambda_\sigma-\Lambda_{\sigma^\prime}\|$.
We use instead the CGO sinogram data-fidelity term
\begin{equation}\label{def:newdatadiscrepancy}
\|\mathcal{S}_\sigma(\theta,\varphi,\rho)-\mathcal{S}_{\sigma^\prime}(\theta,\varphi,\rho)\|_{L^2(\mathbb{T}^2)}^2,
\end{equation}
where $\mathbb{T}^2$ denotes the two-dimensional torus.


%

\subsection{Contrast enhancement}

Assume that the piecewise constant conductivity satisfies $\text{supp}(\sigma-1)\subset\Omega$ and that we know {\em a priori} approximate bounds $0<c<1$ and $C>1$ such that
\begin{eqnarray}\label{eq-apriori_info}
\min_{z\in \Omega}\sigma(z)>c,\qquad \max_{z\in \Omega}\sigma(z)<C.
\end{eqnarray}
Let $\widetilde{\sigma}$ denote an approximate reconstruction to the true conductivity $\sigma$ defined on $\Omega$, whose contrast we intend to improve, and suppose
$\text{supp}(\widetilde{\sigma}-1)\subset\Omega$. Set $f(z)=\widetilde{\sigma}-1$ and denote
\begin{equation}\label{eq-mMdef}
  m=\min_{z\in \Omega}f(z),\qquad M=\max_{z\in \Omega}f(z),
\end{equation}
and assume that $m<0$ and $M>0$. Let $s$ and $t$ be two parameters such that $0\leq s\leq
1$ and $0\leq t\leq 1$ and define
\begin{equation}\label{eq-sigmast}
  \sigma_{s,t}(z) := 1+
\left\{
\begin{array}{ll}
t(C-1)f(z)/M & \mbox{ for } z \mbox{ satisfying } \widetilde{\sigma}(z)>1,\\
s(c-1)f(z)/m & \mbox{ for } z \mbox{ satisfying } \widetilde{\sigma}(z)<1,\\
0 & \mbox{ otherwise.}
\end{array}
\right.
\end{equation}
Note that $\sigma_{0,0}\equiv 1$, and the maximum values $s=1$ and $t=1$
yield the maximal image contrast below and above 1, respectively, i.e.
\begin{eqnarray*}
\min_{z\in \Omega}\sigma_{1,t}(z)=c,\qquad \max_{z\in \Omega}\sigma_{s,1}(z)=C.
\end{eqnarray*}
We determine the optimal values for $s,t$ as the minimizers of the nonlinear data-discrepancy functional based on the CGO sinogram:
\begin{equation}\label{form:minim}
  (s_0,t_0):= \arg\min_{\hspace{-1.5em}\substack{\\(s,t)\in [0,1]^2}}  \frac{\|\mathcal{S}_{\sigma_{s,t}}(\cdot\,,\cdot\,,\rho)-\mathcal{S}_\sigma^{\delta}(\cdot\,,\cdot\,,\rho)\|_{L^2(\mathbb{T}^2)}}{ \|\mathcal{S}_\sigma^{\delta}(\cdot\,,\cdot\,,\rho)\|_{L^2(\mathbb{T}^2)}}.
\end{equation}
The  result of contrast-enhancement is then $\sigma_{CE}:=\sigma_{s_0,t_0}$.

In this paper, the objective function in \eqref{form:minim} is minimized via the DIRECT algorithm \cite{Perttunen}, in an analogous fashion to \cite{Hamilton2014}.  The expression DIRECT refers to ``DIviding RECTangles'', which suggests the strategy of this sampling, global search algorithm. 


\section{Numerical Implementation}\label{sec:numresults}
\subsection{Simulation of noisy EIT data}\label{sec:EITdata}

\begin{figure}
\centering
\begin{picture}(350,135)
\put(20,5){\includegraphics[width=130 pt]{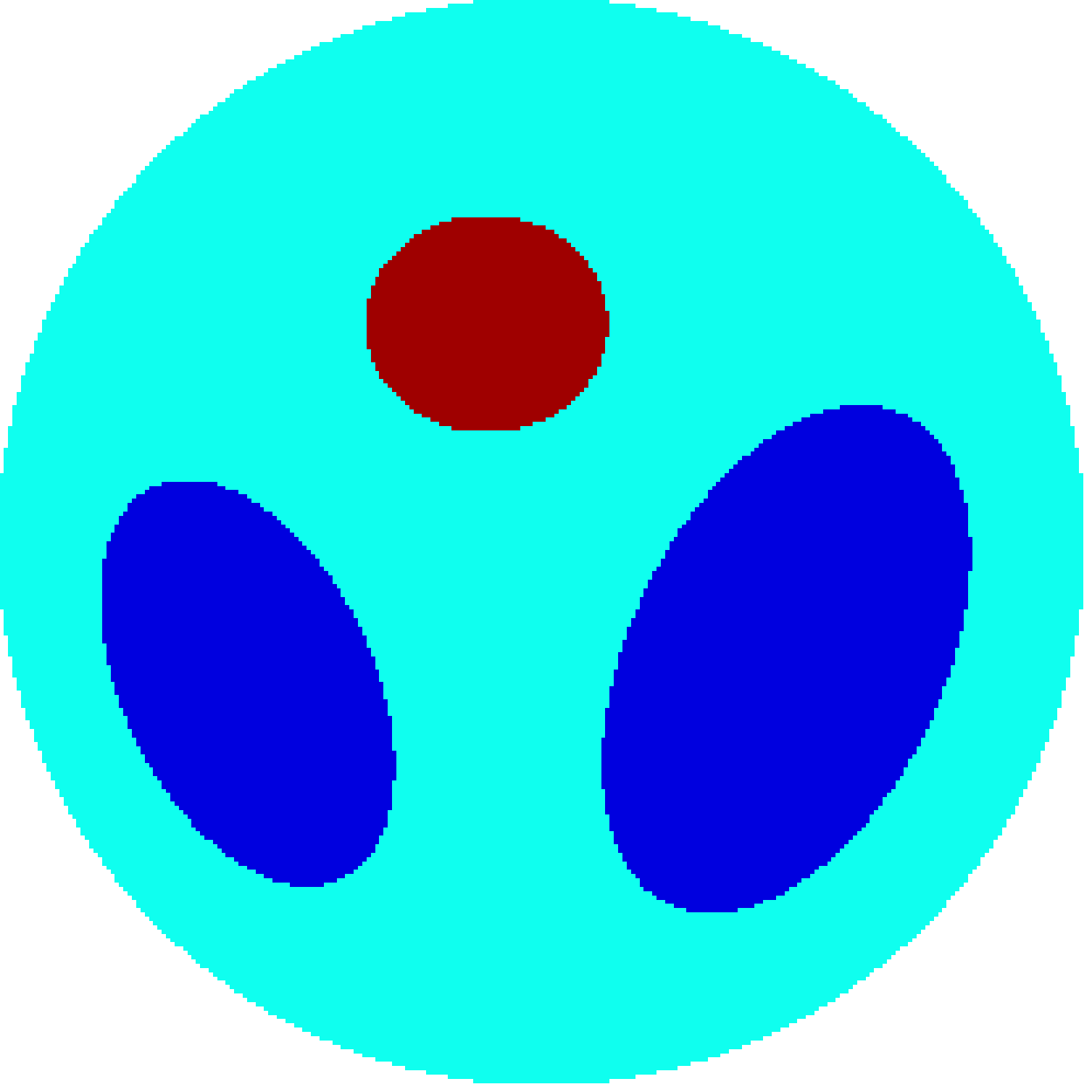}}
\put(170,5){\includegraphics[width=130 pt]{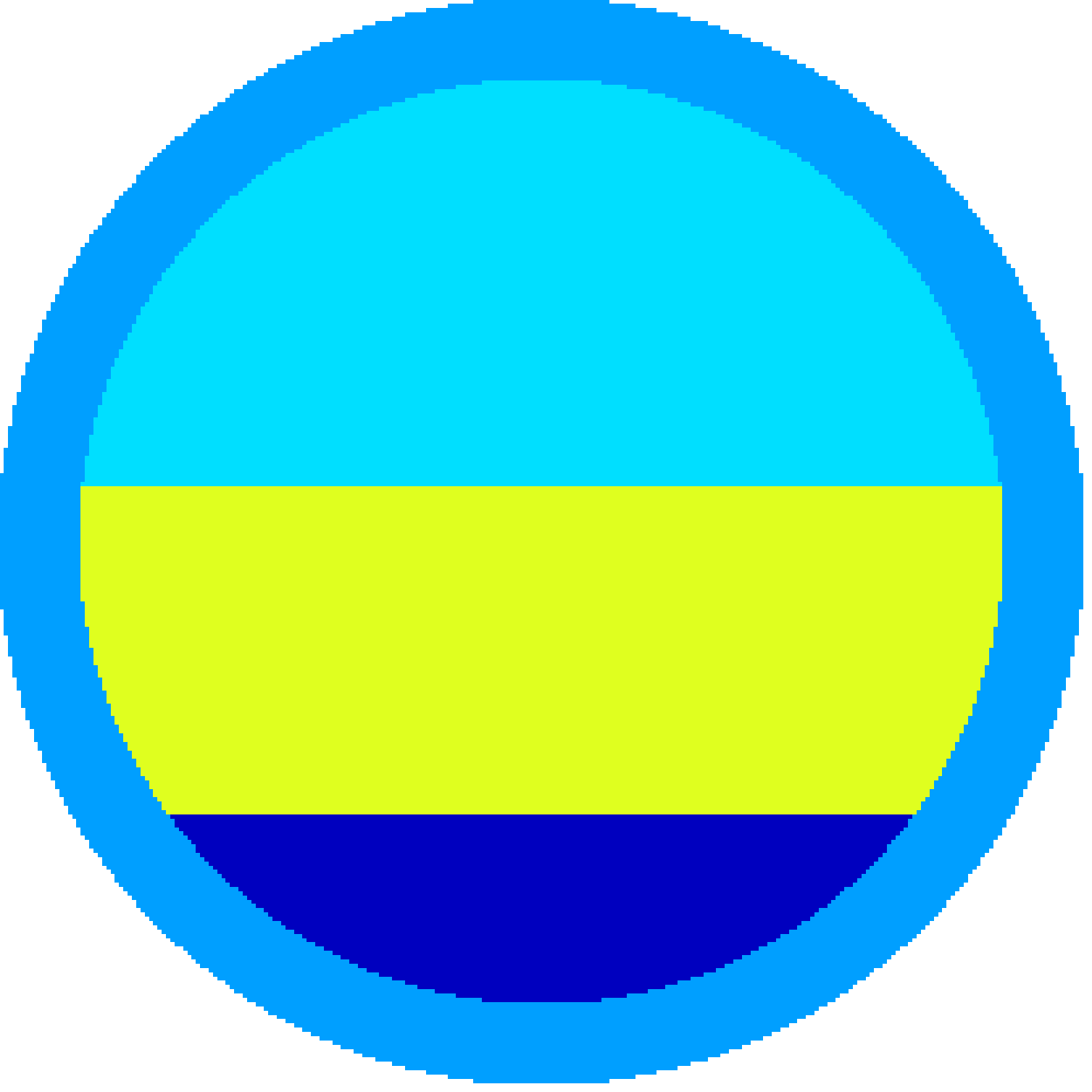}}
\end{picture}
\caption{\label{fig:trueconduc}True conductivity phantoms; for conductivity values see Table \ref{table:conductivities}. Left:
heart-and-lungs phantom $\sigma_1$. Right: pipeline phantom $\sigma_2$.}
\end{figure}

Our numerical experiments deal with two simulated discontinuous conductivity phantoms, namely a heart-and-lungs
phantom $\sigma_1$ and a cross-section of a stratified oil pipeline phantom $\sigma_2$. See Table \ref{table:conductivities} and Figure \ref{fig:trueconduc}.

\begin{table}
\begin{center}
\begin{tabular}{|ll|ll|}
\hline
  $\sigma_1$ & &   $\sigma_2$ & \\
  \hline
  Background & $1.0$ & Pipe & $1.0$ \\
  Lung  & $0.5$ & Top layer (oil) & $1.2$\\
  Heart & $2.0$ & Middle layer (water) & $2.0$ \\
 &&Bottom layer (sand) & $0.3$\\
\hline
\end{tabular}\end{center}
\medskip
\caption{\label{table:conductivities}Conductivity values in the two simulated phantoms shown in  Figure \ref{fig:trueconduc}. Left: heart-and-lungs phantom $\sigma_1$. Right: oil pipeline phantom $\sigma_2$. }
\end{table}

\subsubsection{Computation of the discrete Neumann-to-Dirichlet map}\label{sec:ND}

EIT data was simulated using the finite element method and following \cite[Sections 13.2.3 and 16.3.3]{Mueller2012}.  We use the trigonometric basis functions
\[\phi_n(\theta) = \left\{
\begin{array}{ll}
\pi^{-1/2}\cos((n+1)\theta/2) &\mbox{for odd }n,\\
\pi^{-1/2}\sin(n\theta/2) &\mbox{for even }n,\\
\end{array}
\right.\]
where $1\leq n\leq 2N$. The
Neumann-to-Dirichlet map $\mathcal{R}_{\sigma}$ is approximated by
the matrix $\textbf{R}_{\sigma}=[(\textbf{R}_{\sigma})_{m,n}]$
given by
\[(\textbf{R}_{\sigma})_{m,n} = \langle \textbf{R}_{\sigma}\phi_n ,
\phi_m \rangle = \int_0^{2\pi}
(\textbf{R}_{\sigma}\phi_n)(\theta)\,\phi_m(\theta)\, d\theta,\]
where $\textbf{R}_{\sigma}\phi_n = u_n|_{\partial\Omega}$ with
$\nabla\cdot \sigma\nabla u_n =0$ in $\Omega$, $(\sigma\,(\partial
u_n / \partial\nu))|_{\partial\Omega} = \phi_n$ and
$\int_{\partial\Omega}u_n\, dS = 0$. Here, $1\leq m,n\leq 2N$.  In this paper, we use $N=16$ which corresponds to 33 linearly independent current patterns.


Relative Gaussian noise was added to the boundary
voltage data as in \cite{Hamilton2014a}. Namely,
define $\widetilde{\textbf{R}}_{\sigma}$ with
$(\widetilde{\textbf{R}}_{\sigma})_{m,n} = \langle
\widetilde{\textbf{R}}_{\sigma}\,\phi_n , \phi_m \rangle$, where
\begin{equation}\label{addnoise}
\widetilde{\textbf{R}}_{\sigma}\,\phi_n =
\textbf{R}_{\sigma}\phi_n + \eta\,
\mathcal{N}_n\norm{\textbf{R}_{\sigma}\phi_n}_{L^{\infty}},
\end{equation}
$\eta$ denotes the noise level so that $100\eta\%$ noise is added,
and $\mathcal{N}_1$,...,$\mathcal{N}_{2N+1}$  are independent
Gaussian distributions with mean zero and variance one, which are
implemented through the MATLAB function {\it randn} generating
pseudorandom values drawn from the standard normal distribution
$\mathcal{N}(0,1)$.  The noisy D-N matrix is then formed by inverting
$\widetilde{\textbf{R}}_{\sigma}$ and adding a zero row and column on top and left as described in \cite{Mueller2012}.

The algorithm outlined in Figure \ref{fig:flowchart} was
applied to the two  phantoms shown in Figure~\ref{fig:trueconduc}.  In this paper, the algorithm was performed for $J=3$ iterations, and for three noise levels as follows: zero added noise, $0.1\%$ added noise and $0.75\%$ added noise (corresponding to setting  $\eta=0,\,0.001,\,0.0075$).  The output of the algorithm after the $J=3$ iterations is denoted by $\sigma_{\mbox{\tiny CE}}^{(J)}$.


\subsection{Computational grids}\label{sec:grids}

All the computations on the $z$-plane (D-bar reconstructions, TV flow outputs and CE outputs) were generated on a $z$-grid of $2^{\ell}{\times}2^{\ell}$ equidistributed points of the square $[-s,s)\times [-s,s)$ with $\ell=8$ and $s=2.3$.
Thus, the $z$-grid consists of $2^{16}=65536$ points. The D-bar equation solver (used to solve \eqref{eq-dbark-mu-sol}) was executed just on the set of $9729$ points of the $z$-grid belonging to the closed disc $\overline{\Omega}$.  Note that the larger $z$-region is needed to extend the scattering data.

The $k$-grids are problem specific, i.e. for lower levels of noise a larger radius $R$ can be used for the initial low-pass filtering in the nonlinear Fourier domain.  In each case, we fixed the two parameters $R$ (the initial low pass filtering chosen intuitively by looking at where the scattering data ``blows up'' in magnitude), and $\widetilde{R}$ (the increased scattering radius to be determined by solving the Beltrami equation).  In each case, the $k$-grid for the scattering data was comprised of $2^7\by 2^7$ equispaced points on the square $[-\widetilde{R},\widetilde{R})\by [-\widetilde{R},\widetilde{R})$.  The traditional (and original) D-bar image $\sigma_{\mbox{\tiny DB}}^{(1)}$ is computed from scattering data satisfying $|k|\leq R$ and all subsequent D-bar images $\sigma_{\mbox{\tiny DB}}^{(2)}$, $\sigma_{\mbox{\tiny DB}}^{(3)}$, etc. are computed using the larger disc $|k|\leq \widetilde{R}$.

\vspace{1em}
\noindent{\it Remark}:  If supp$(\sigma-\sigma_0)\subset D(0,1)$ for some positive constant $\sigma_0\neq 1$, the algorithm can be re-scaled as follows. Defining
$\widetilde{\sigma}:=\sigma/\sigma_0$, we have
supp$(\widetilde{\sigma}-1)\subset D(0,1)$. Apply the above
algorithm to $\Lambda_{\widetilde{\sigma}}=\sigma_0
\Lambda_{\sigma}$ and write $\widetilde{\sigma}_{\mbox{\tiny
CE}}^{(J)}$ for the output generated. Take
$\sigma_0\,\widetilde{\sigma}_{\mbox{\tiny CE}}^{(J)}$ as the
final approximation to $\sigma$.

\subsection{Computation of the initial nonlinear Fourier transform}

This step corresponds to Figure \ref{fig:flowchart}(a).
From the boundary measurements $\Lambda_\sigma^\delta$ we use \eqref{eq:CGOBeltrami_traces_Mmu} and
\eqref{eq:tau_Beltrami_from_Mu} to
compute the initial scattering data $\tau^0(k)$ on
a $k$-disc of radius $R$. We write $\T^0(k):= -4\pi i \bar{k} \tau^0(k)$.

\subsection{D-bar reconstruction}

This step corresponds to Figure \ref{fig:flowchart}(b). Use the shortcut method described in Section \ref{sec:shortcutmethod} to compute $\sigma^{(j)}_{\mbox{\tiny DB}}$ from the scattering transform $\T^{j-1}(k)$. If $j=1$ then use the smaller cutoff disc $|k|<R$, else use larger disc $|k|<\widetilde{R}$.

\subsection{Edge-enhancement using TV flow}

This step corresponds to Figure \ref{fig:flowchart}(c). Introduce edges into the D-bar image
$\sigma^{(j)}_{\mbox{\tiny DB}}$ by applying the segmentation flow of Section \ref{sec:TVflow}. The resulting piecewise constant image, defined by \eqref{eq:segregularized}, is called $\sigma_{\mbox{\tiny TV}}^{(j)}$. Note that the initial guess of the mean intensity are obtained directly by K-means algorithm, where $K$ denote the number regions pre-selected. In practice, this is a reasonable guess at the number of regions of different conductivity in your domain.  

\subsection{Contrast enhancement}

This step corresponds to Figure \ref{fig:flowchart}(d). The CGO sinograms were implemented via $33\times 33$ matrices
$[(\textbf{S})_{m,l}]$ with $(\textbf{S})_{m,l} =
M_{\mu}(z_m,k_l)-1$ $=M_{\mu}(e^{i\theta_m},2 e^{i\varphi_l})-1$,
where $$\theta_m = (m-1-N)2\pi/(2N+1),\qquad 1\leq m\leq
2N+1,$$ with $N=16$, $\varphi_l = \theta_l$, and $M_{\mu}(z,k)$ refers to the
solution explained in Section \ref{sec:CGOsino} for both the true
$\sigma$ and the corresponding approximations $\sigma_{s,t}$. Therefore, $\{\theta_m\}$ and $\{\varphi_l\}$ are the same partition of the interval $(-\pi,\pi)$.  Note that one can choose the points $\varphi_l$ independently of the $\theta$ values if desired.

Finally, the output image $\sigma_{\mbox{\tiny CE}}^{(j)}$ is determined by
plugging $\widetilde{\sigma}=\sigma^{(j)}_{\mbox{\tiny TV}}$ into
\eqref{eq-sigmast} for $(s_0,t_0)\in [0,1]^2$ obtained via the
DIRECT optimization strategy (see \eqref{form:minim}).

\subsection{Extension of the Scattering Transform}\label{sec-extendScat}

This step corresponds to Figure \ref{fig:flowchart}(e). The radius of the admissible scattering data is increased from $R$ to $\tilde{R}$ by computing new stable scattering data in the annulus $R-1<|k|<\tilde{R}$ corresponding to the TV-sharpened and contrast adjusted image $\sigma_{\mbox{\tiny CE}}^{(j)}$ as follows.  First, evaluate the Beltrami coefficient $\mu(z)=(1-\sigma_{\mbox{\tiny
CE}}^{(j)}(z))/(1+\sigma_{\mbox{\tiny CE}}^{(j)}(z))$ on the $z$-grid $[-2.3,2.3)\by[-2.3,2.3)$.  Next, solve the Beltrami equation \eqref{eq:BeltramiEq} for the {\sc CGO} solutions $f_{\pm\mu}(z,k)$ for $k$ in $k$-annulus $R-1\leq |k|\leq \widetilde{R}$.  Finally, evaluate the scattering data $\tau(k)$ in the $k$-annulus $R-1\leq |k|\leq \widetilde{R}$ via
\begin{equation}\label{eq:scat_AP}
\overline{\tau(k)} := {1\over 2\pi} \int_{\R^2}
\left(\overline{\partial}_z \middle(M_{+\mu}(z,k) - M_{-\mu}(z,k) \right)
dz_1 dz_2\,
\end{equation}
where $M_{\pm\mu}(z,k)=e^{-ikz}f_{\pm\mu}(z,k)$.  Call this new scattering data $\tilde{\tau}^{(j)}(k)$.
The new scattering data on the larger radius is then
\[\tau^{(j)}(k):=\chi(k)\tau^{0}(k)+\left(1-\chi(k)\right)\tilde{\tau}^{(j)}(k), \quad |k|<\tilde{R},\]
where we used the polynomial radial cutoff function $\chi$ defined by
\[\chi(k) = \left\{
\begin{array}{ll}
1, & \mbox{if }|k|<R-1,\\
p(|k|-(R-1)), & \mbox{if }R-1<|k|<R,\\
0, & \mbox{if }|k|>R,\\
\end{array}
\right.\]
with $p(t) = 1-3 t^2 + 2 t^3$ to blend the data in the overlap region $R-1\leq |k|\leq R$. Note that in Figure \ref{fig:flowchart} we use the notation $\T^j(k):= -4\pi i \bar{k} \tau^{(j)}(k)$ and $\tilde{\T}^j(k):= -4\pi i \bar{k} \tilde{\tau}^{(j)}(k)$.

\subsubsection{Computation of ``True'' Scattering Data}
We need a comparison for our new extended scattering data with the best possible scattering data.  By \textit{best possible scattering data} we refer to the scattering data that is obtained by computing the CGO solutions to the Beltrami equation \eqref{eq:BeltramiEq} with $\mu$ corresponding to the true $\sigma$, and evaluating the scattering data $\tau(k)$ via \eqref{eq:scat_AP}. For details on how to solve the Beltrami equation and generate the scattering data $\tau$, the reader is referred to \cite{Astala2010,Astala2014,Huhtanen2012}.

\section{Numerical Results}\label{sec:numresults2}
The algorithm was tested on the two phantoms shown in Figure~\ref{fig:trueconduc} and the results are shown here.

\subsection{Example 1: A Heart and Lungs Phantom}
For the heart and lungs phantom, $\sigma_1$, it was assumed known {\em apriori} that the internal conductivity was bounded between $c=0.3$ and $c=2.5$.  Setting the initial scattering radius $R$ was 5 and the enlarged radius $\widetilde{R}$ to 10 for all noise levels proved sufficient.  The initial scattering data was reliable for all noise levels within the $k$-disc of radius 5.  The parameters for the TV flow were $K=4$ and $\lambda=0.1$.

\begin{figure}[h!]
\centering
\begin{picture}(350,350)
\put(0,10){\includegraphics[width=330pt]{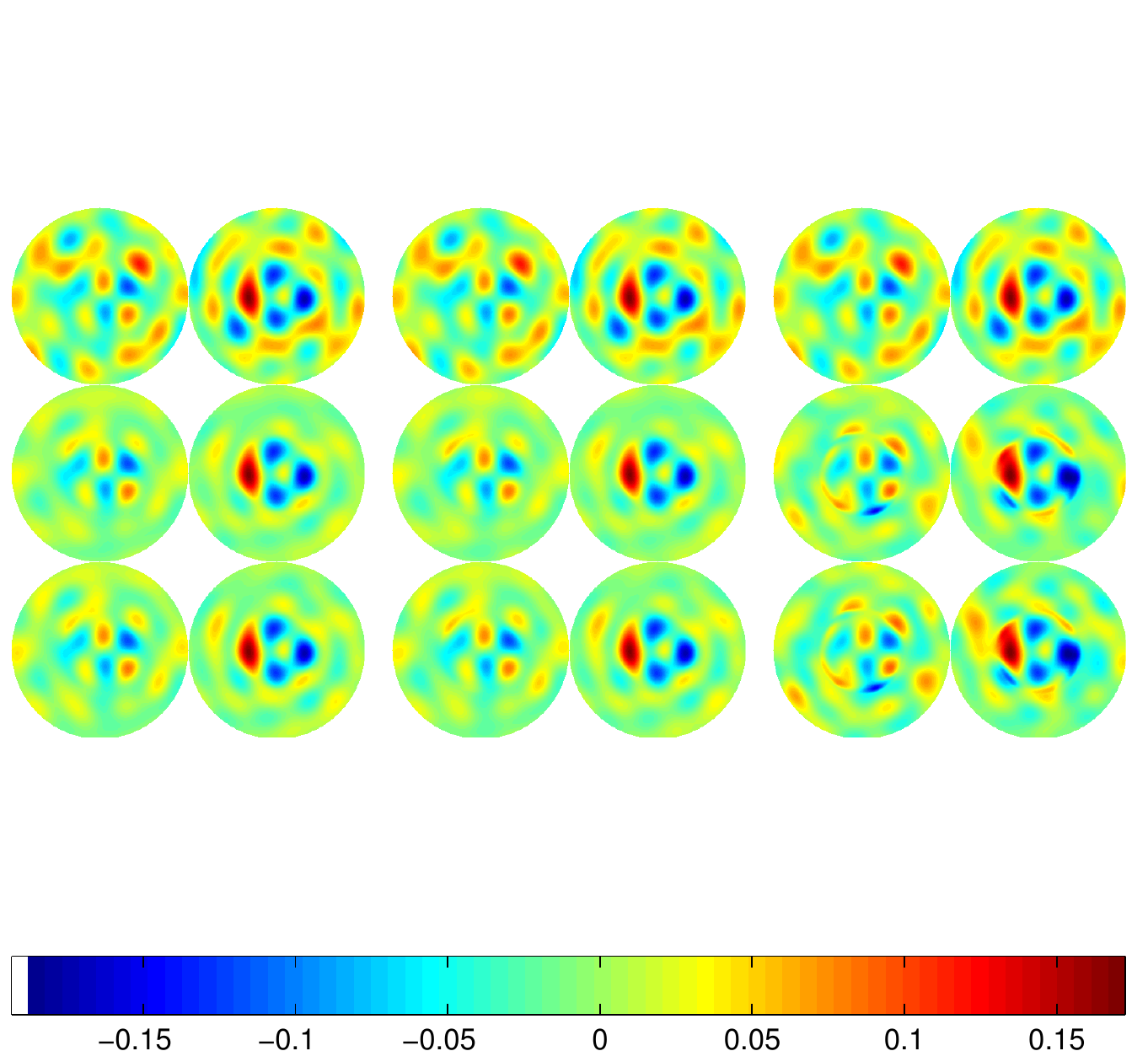}}
\put(-20,120){$\tau^{(2)}$}
\put(-20,170){$\tau^{(1)}$}
\put(-20,220){$\tau_{\tiny \mbox{B}}$}
\put(20,263){\textbf{Re}}
\put(70,263){\textbf{Im}}
\put(130,263){\textbf{Re}}
\put(180,263){\textbf{Im}}
\put(240,263){\textbf{Re}}
\put(290,263){\textbf{Im}}
\put(10,75){a) no added noise}
\put(125,75){b) $0.1\%$ noise}
\put(235,75){c) $0.75\%$ noise}
\end{picture}
\caption{\label{fig:allscat_hl}Images of the real and imaginary parts of the reliable scattering transform $\tau_{\mbox{\tiny B}}$ (computed directly from the Beltrami equation) and the combined scattering data for the two iterations $\tau^{(1)}$ and $\tau^{(2)}$ from the simulated Dirichlet-to-Neumann corresponding to the heart and lungs phantom $\sigma_1$.}
\end{figure}
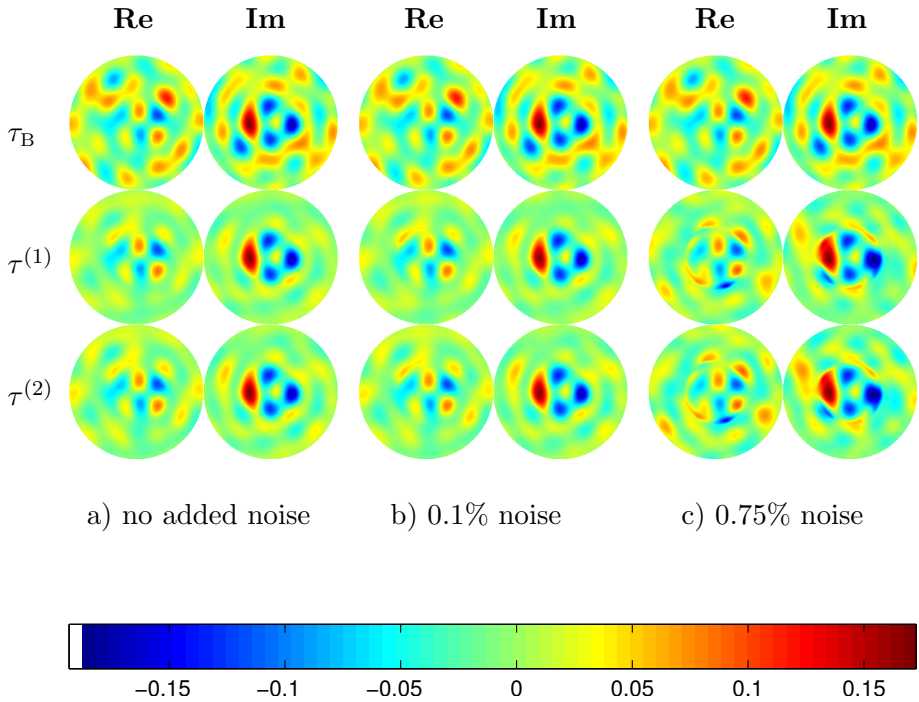

The scattering data for each noise level is displayed in Figure~\ref{fig:allscat_hl}.  The figures contain images of the actual Beltrami scattering transform $\tau_{\mbox{\tiny B}}$ on the larger $k$ disc of radius 10.  This is used to evaluate the efficacy of the new proposed approach.  The \textit{true} scattering data $\tau_{\mbox{\tiny B}}$ was computed by solving \eqref{eq:BeltramiEq} and \eqref{eq:scat_AP} with the known $\mu=\frac{1-\sigma_1}{1+\sigma_1}$.  It is to serve as a \textit{best case scenario} baseline.

\begin{figure}[h!]
\centering
\begin{picture}(350,270)
\put(135,200){\includegraphics[width=50pt]{figure12.pdf}}
\put(135,255){\textbf{True $\sigma_1$}}
\put(0,40){\includegraphics[width=110pt]{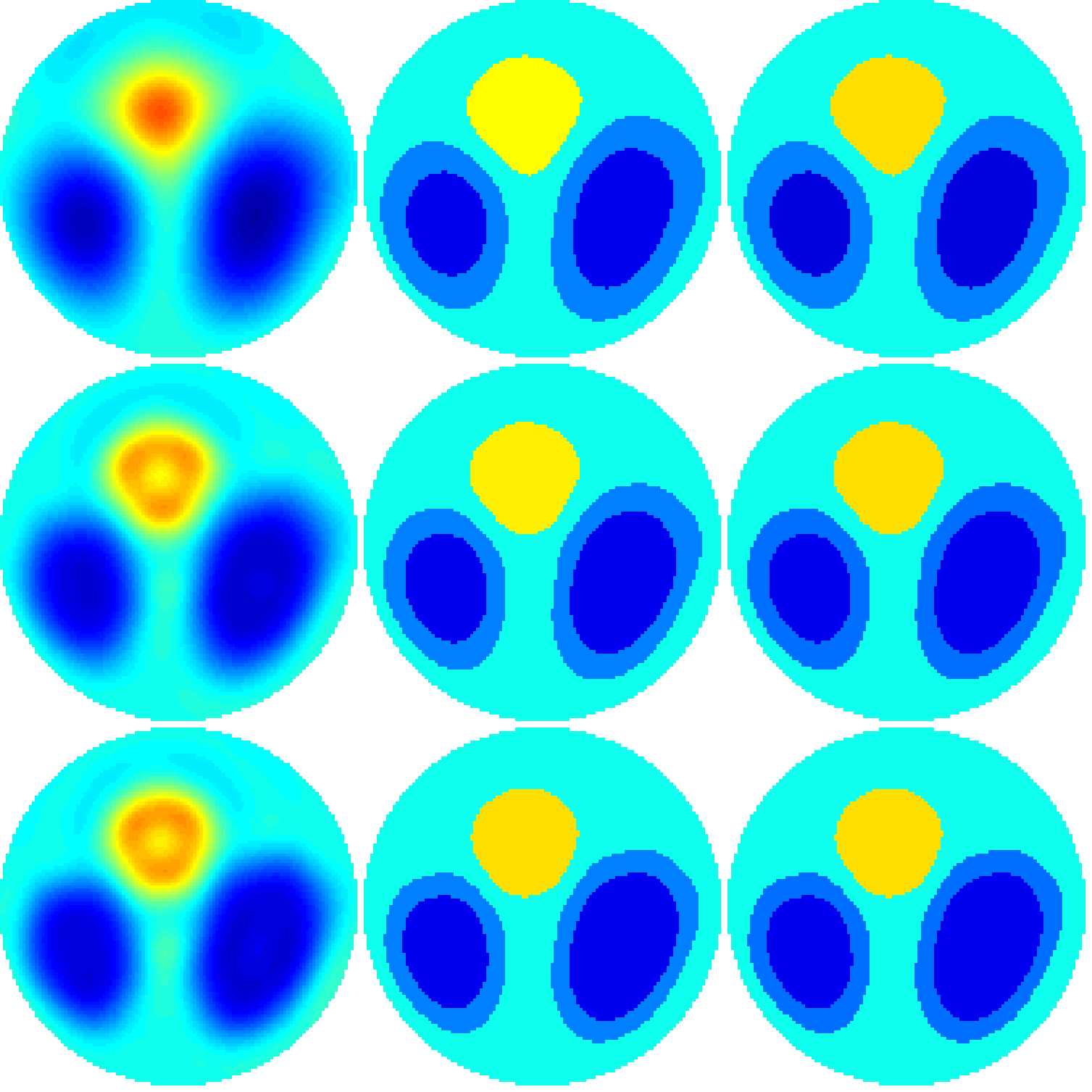}}
\put(236,38){\includegraphics[width=145pt]{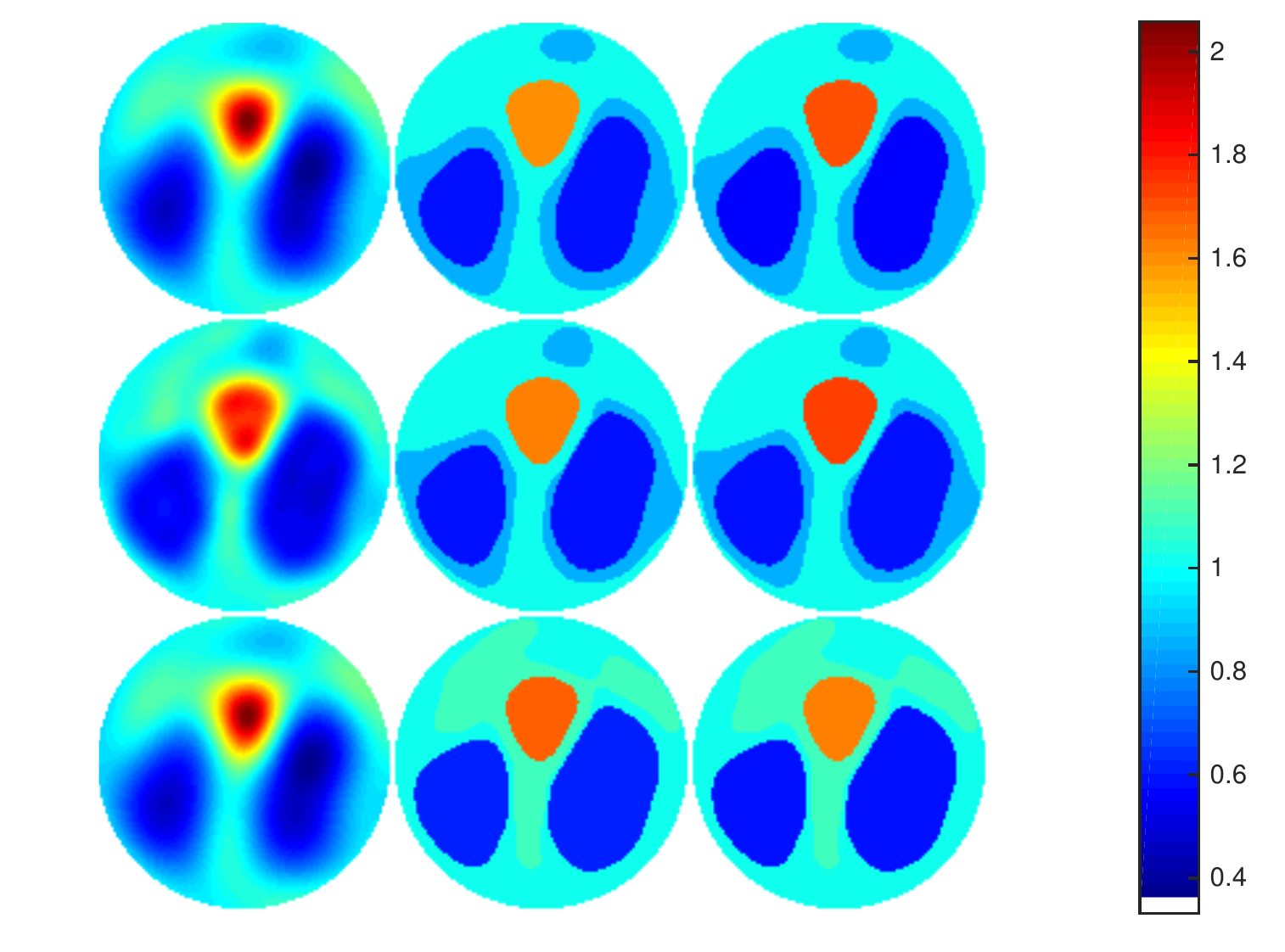}}
\put(117,40){\includegraphics[width=110pt]{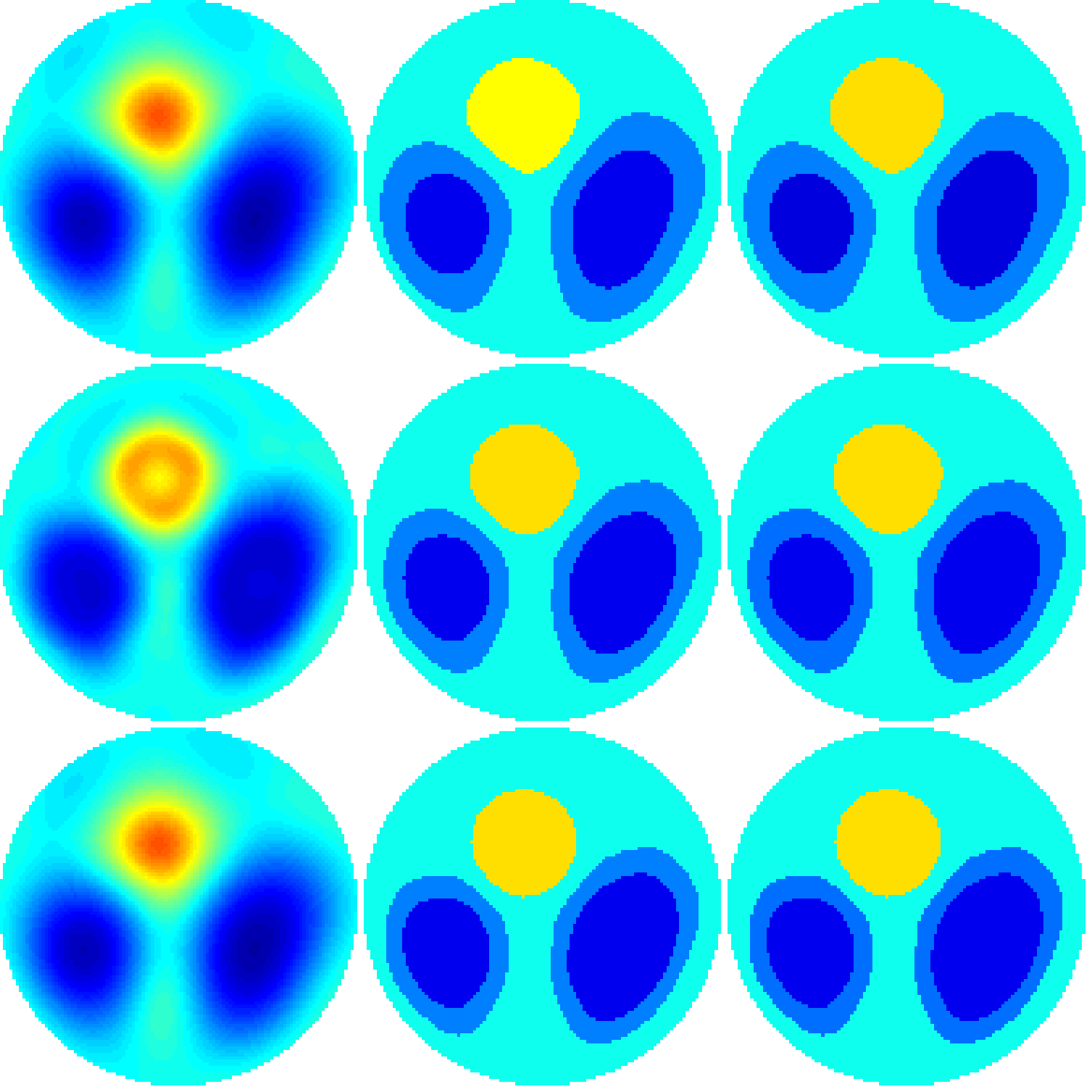}}

\put(-40,130){Iter 1}
\put(-40,95){Iter 2}
\put(-40,58){Iter 3}
\put(6,160){$\sigma_{\tiny \mbox{DB}}$}
\put(41,160){$\sigma_{\tiny \mbox{TV}}$}
\put(80,160){$\sigma_{\tiny \mbox{CE}}$}

\put(125,160){$\sigma_{\tiny \mbox{DB}}$}
\put(160,160){$\sigma_{\tiny \mbox{TV}}$}
\put(200,160){$\sigma_{\tiny \mbox{CE}}$}

\put(245,160){$\sigma_{\tiny \mbox{DB}}$}
\put(281,160){$\sigma_{\tiny \mbox{TV}}$}
\put(320,160){$\sigma_{\tiny \mbox{CE}}$}

\put(0,5){a) no added noise}
\put(130,5){b) $0.1\%$ noise}
\put(250,5){c) $0.75\%$ noise}
\end{picture}
\caption{\label{fig:all3recon_hl}This figure shows the real parts of the numerical approximations obtained for the heart and lungs phantom of Example~1, i.e. $\sigma_1$.  The picture consists of three parts a), b), c) corresponding to our three cases of added simulated noise in the EIT voltage data: zero added noise, noise of relative amplitude $0.1\%$ and noise of relative amplitude $0.75\%$, respectively. For the sake of comparison, the true conductivity $\sigma_1$ is displayed above the reconstructions (all on the same color scale).  For each noise level, the D-bar reconstruction $\sigma^{(j)}_{\mbox{\tiny DB}}$ (left columns), the TV sharpened image $\sigma^{(j)}_{\mbox{\tiny TV}}$ (middle columns), and the contrast adjusted TV sharpened images $\sigma^{(j)}_{\mbox{\tiny CE}}$. The first, second, and third rows correspond to the first, second and third iterations, respectively.}
\end{figure}

The reconstructed conductivities for each stage of the algorithm are displayed in Figure~\ref{fig:all3recon_hl}.  Note that the reconstructions are displayed on the same color scale as the original conductivity shown in Figure~\ref{fig:trueconduc} (left) for ease of comparison.

Tables~\ref{table:Errors_0o0noise_HnL}, \ref{table:Errors_0o1noise_HnL}, and \ref{table:Errors_0o75noise_HnL} show the relative $L^2$ errors and the Structural SiMilarity Index (SSIM) values for the heart-and-lungs phantom $\sigma_1$ of Example~1 for zero added noise, $0.1\%$ added noise, and $0.75\%$ added noise, respectively.  The error values are presented for each step of the proposed algorithm: the D-bar reconstruction $\sigma^{(j)}_{\mbox{\tiny DB}}$, the sharpened reconstruction $\sigma^{(j)}_{\mbox{\tiny TV}}$, and the contrast adjusted sharpened reconstruction $\sigma^{(j)}_{\mbox{\tiny CE}}$.

\begin{table}[h!]
\centering
{\footnotesize
\begin{tabular}{l|lllll|lll}

  \multicolumn{4}{c}{$L^2$ Relative Error} & $\quad$ & \multicolumn{4}{c}{SSIM} \\
  \multicolumn{9}{c}{}\\

  $j$ & $\sigma^{(j)}_{\mbox{\tiny DB}}$ & $\sigma^{(j)}_{\mbox{\tiny TV}}$ & $\sigma^{(j)}_{\mbox{\tiny CE}}$
   & & $j$ & $\sigma^{(j)}_{\mbox{\tiny DB}}$ & $\sigma^{(j)}_{\mbox{\tiny TV}}$ & $\sigma^{(j)}_{\mbox{\tiny CE}}$ \\

  \cline{1-4} \cline{6-9}

1 & 0.1240 & 0.1240  & 0.1202 & &
1 & 0.6600 & 0.6600 & 0.6541 \\

2 & 0.1095 & 0.1185 & 0.1168 & &
2 & 0.7351 & 0.6903 & 0.6878 \\

 3 & 0.1054 & 0.1157 & 0.1145 & &
3 & 0.7425& 0.7117 & 0.7096 \\

\end{tabular}
}
\caption{\label{table:Errors_0o0noise_HnL}The $L^2$ relative errors as well as the SSIM values for the zero added noise case for Example~1 with the heart and lungs phantom $\sigma_1$.}
\end{table}

\begin{table}[h!]
\centering
{\footnotesize
\begin{tabular}{l|lllll|lll}

  \multicolumn{4}{c}{$L^2$ Relative Error} & $\quad$ & \multicolumn{4}{c}{SSIM} \\
  \multicolumn{9}{c}{}\\

  $j$ & $\sigma^{(j)}_{\mbox{\tiny DB}}$ & $\sigma^{(j)}_{\mbox{\tiny TV}}$ & $\sigma^{(j)}_{\mbox{\tiny CE}}$
   & & $j$ & $\sigma^{(j)}_{\mbox{\tiny DB}}$ & $\sigma^{(j)}_{\mbox{\tiny TV}}$ & $\sigma^{(j)}_{\mbox{\tiny CE}}$ \\

  \cline{1-4} \cline{6-9}

1 & 0.1009 & 0.1233  & 0.1194 & &
1 & 0.7304 & 0.6603 & 0.6545 \\

2 & 0.1083 & 0.1174 & 0.1158 & &
2 & 0.7348 & 0.6907 & 0.6883 \\

 3 & 0.1009 & 0.1142 & 0.1133 & &
3 & 0.7304& 0.7131 & 0.7111 \\

\end{tabular}
}
\caption{\label{table:Errors_0o1noise_HnL}The $L^2$ relative errors as well as the SSIM values for the $0.1\%$ added noise case for Example~1 with the heart and lungs phantom $\sigma_1$.}
\end{table}

\begin{table}[h!]
\centering
{\footnotesize
\begin{tabular}{l|lllll|lll}

  \multicolumn{4}{c}{$L^2$ Relative Error} & $\quad$ & \multicolumn{4}{c}{SSIM} \\
  \multicolumn{9}{c}{}\\

  $j$ & $\sigma^{(j)}_{\mbox{\tiny DB}}$ & $\sigma^{(j)}_{\mbox{\tiny TV}}$ & $\sigma^{(j)}_{\mbox{\tiny CE}}$
   & & $j$ & $\sigma^{(j)}_{\mbox{\tiny DB}}$ & $\sigma^{(j)}_{\mbox{\tiny TV}}$ & $\sigma^{(j)}_{\mbox{\tiny CE}}$ \\

  \cline{1-4} \cline{6-9}

1 & 0.1092 & 0.1202  & 0.1167 & &
1 & 0.6897 & 0.6680 & 0.6639 \\

2 & 0.1076 & 0.1165 & 0.1134 & &
2 & 0.6964 & 0.6816 & 0.6776 \\

 3 & 0.1092 & 0.1151 & 0.1171 & &
3 & 0.6897& 0.6773 & 0.6790 \\

\end{tabular}
}
\caption{\label{table:Errors_0o75noise_HnL}The $L^2$ relative errors as well as the SSIM values for the $0.75\%$ added noise case for Example~1 with the heart and lungs phantom $\sigma_1$.}
\end{table}

In the case of zero and $0.1\%$ added relative noise, the $L^2$ relative error decreases with each iteration and the SSIM value increases.  Both measures confirm that the image is ``improving''.  In the case of $0.75\%$ added relative noise, the $L^2$ relative error and SSIM remain approximately the same but visually one can see the improvements in each iteration as the artifact present in iterations 1 and 2 (above the heart) is absent in iteration 3.

\subsection{Example 2: An Industrial Pipeline}
For the cross-section of the industrial pipeline phantom, $\sigma_2$, it was assumed known {\em apriori} that the internal conductivity was bounded between $c=0.1$ and $c=2.5$. Here we used The parameters for the TV flow were $K=5$ in the TV flow and Table~\ref{table:PipelineParameters} gives the values of the scattering radii and the $\lambda$ parameter used in for the TV flow for each noise level.

\begin{table}[h!]
\centering
\begin{tabular}{|l|l|l|l|}
\hline

{Added Noise Level} & $0\%$ & $0.1\%$ & $0.75\%$ \\
\hline
\hline
$R$               & 6     & 5       & 4        \\
\hline
$\widetilde{R}$   & 10    & 8.3     & 6.6      \\
\hline
$\lambda$         & 0.3   & 0.5     & 0.5     \\
\hline

\end{tabular}
\vspace{0.5em}
\caption{Parameter values for the new algorithm for each noise level for the industrial pipeline phantom $\sigma_2$.}\label{table:PipelineParameters}
\end{table}

\begin{figure}[h!]
\centering
\begin{picture}(350,250)
\put(4,0){\includegraphics[width=330pt]{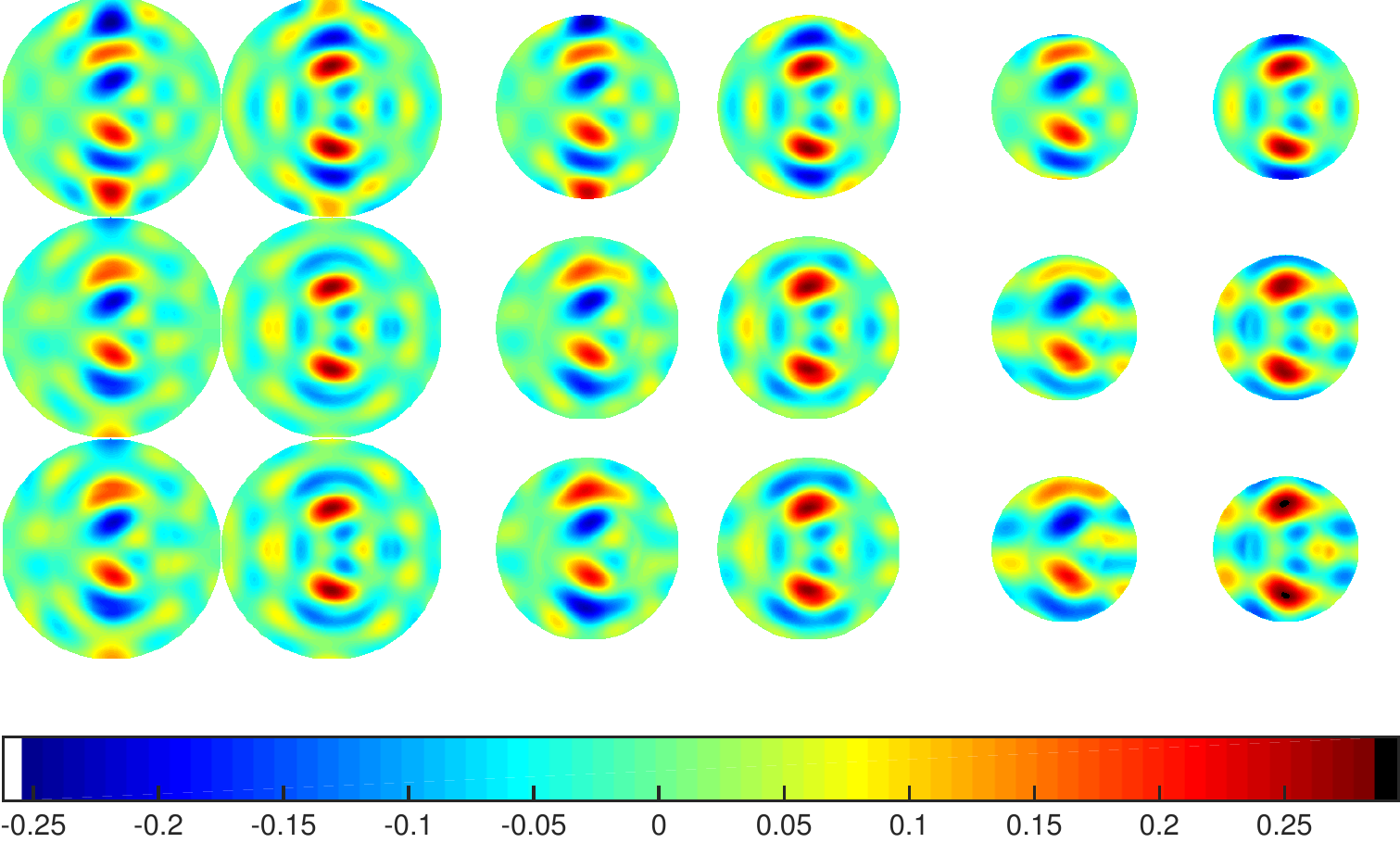}}
\put(-20,70){$\tau^{(2)}$}
\put(-20,120){$\tau^{(1)}$}
\put(-20,170){$\tau_{\tiny \mbox{B}}$}
\put(20,213){\textbf{Re}}
\put(70,213){\textbf{Im}}
\put(130,213){\textbf{Re}}
\put(182,213){\textbf{Im}}
\put(242,213){\textbf{Re}}
\put(294,213){\textbf{Im}}
\end{picture}
\caption{\label{fig:allscat_pipeline}Images of the real and imaginary parts of the reliable scattering transform $\tau_{\mbox{\tiny B}}$ (computed directly from the Beltrami equation) and the combined scattering data for the two iterations $\tau^{(1)}$ and $\tau^{(2)}$ from the simulated Dirichlet-to-Neumann corresponding to the industrial pipeline phantom $\sigma_2$. Left: zero added noise and $\tilde{R}=10$, Middle: $0.1\%$ added noise and $\tilde{R}=8.3$, Right: $0.75\%$ added noise and $\tilde{R}=6.6$.}
\end{figure}
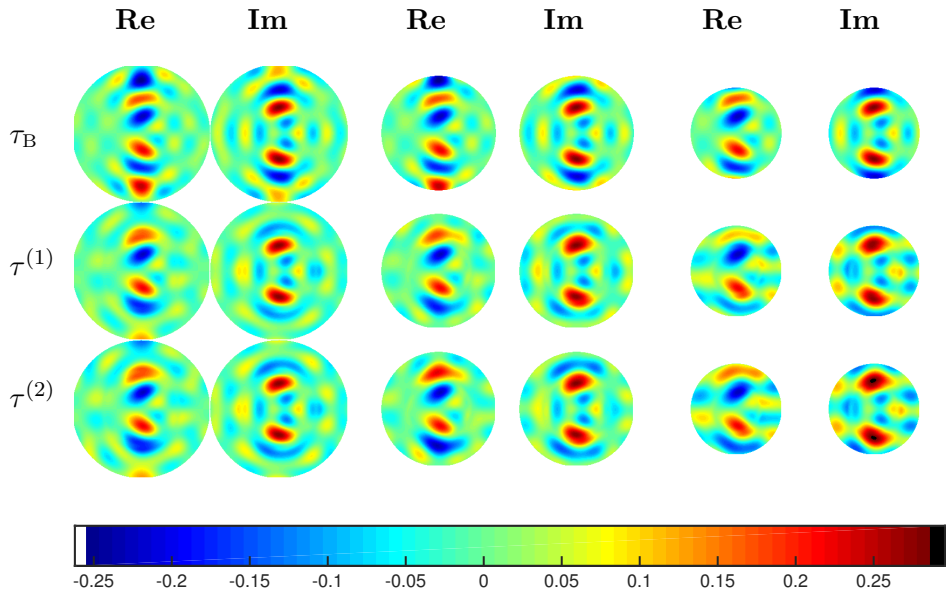

The scattering data for each noise level is displayed in Figure~\ref{fig:allscat_pipeline}. Here the figures contain images of the actual Beltrami scattering transform $\tau_{\mbox{\tiny B}}$ on the larger $k$ disc of radii of $10,\;8.3,\;6.6$ for the varying noise levels respectively.  This is used to evaluate the efficacy of the new proposed approach.  The \textit{true} scattering data $\tau_{\mbox{\tiny B}}$ was computed by solving \eqref{eq:BeltramiEq} and \eqref{eq:scat_AP} with the known $\mu=\frac{1-\sigma_2}{1+\sigma_2}$.  It is to serve as a \textit{best case scenario} baseline.

\begin{figure}[h!]
\centering
\begin{picture}(350,270)

\put(135,200){\includegraphics[width=50pt]{figure13.pdf}}
\put(135,255){\textbf{True $\sigma_2$}}
\put(0,40){\includegraphics[width=110pt]{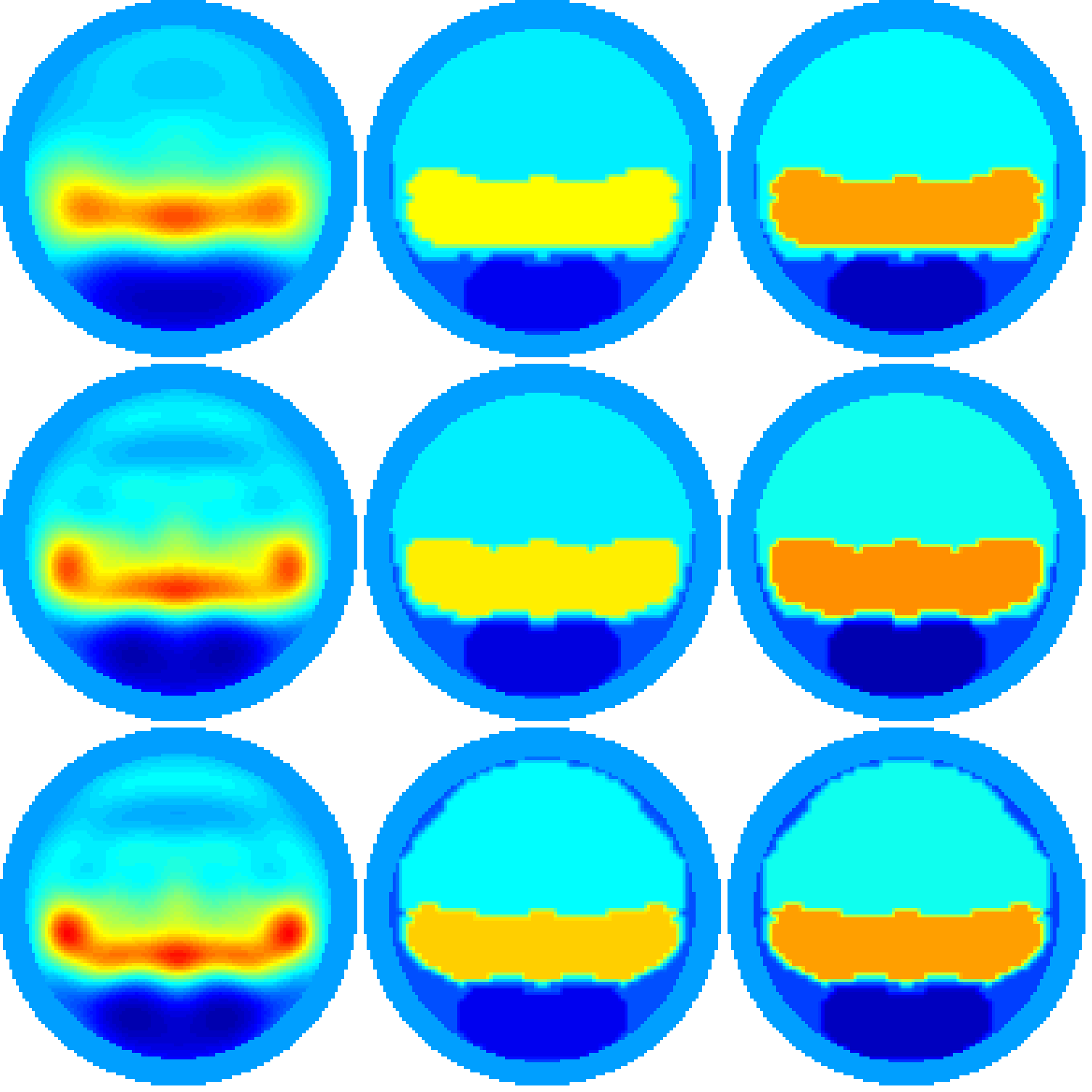}}

\put(236,38){\includegraphics[width=145pt]{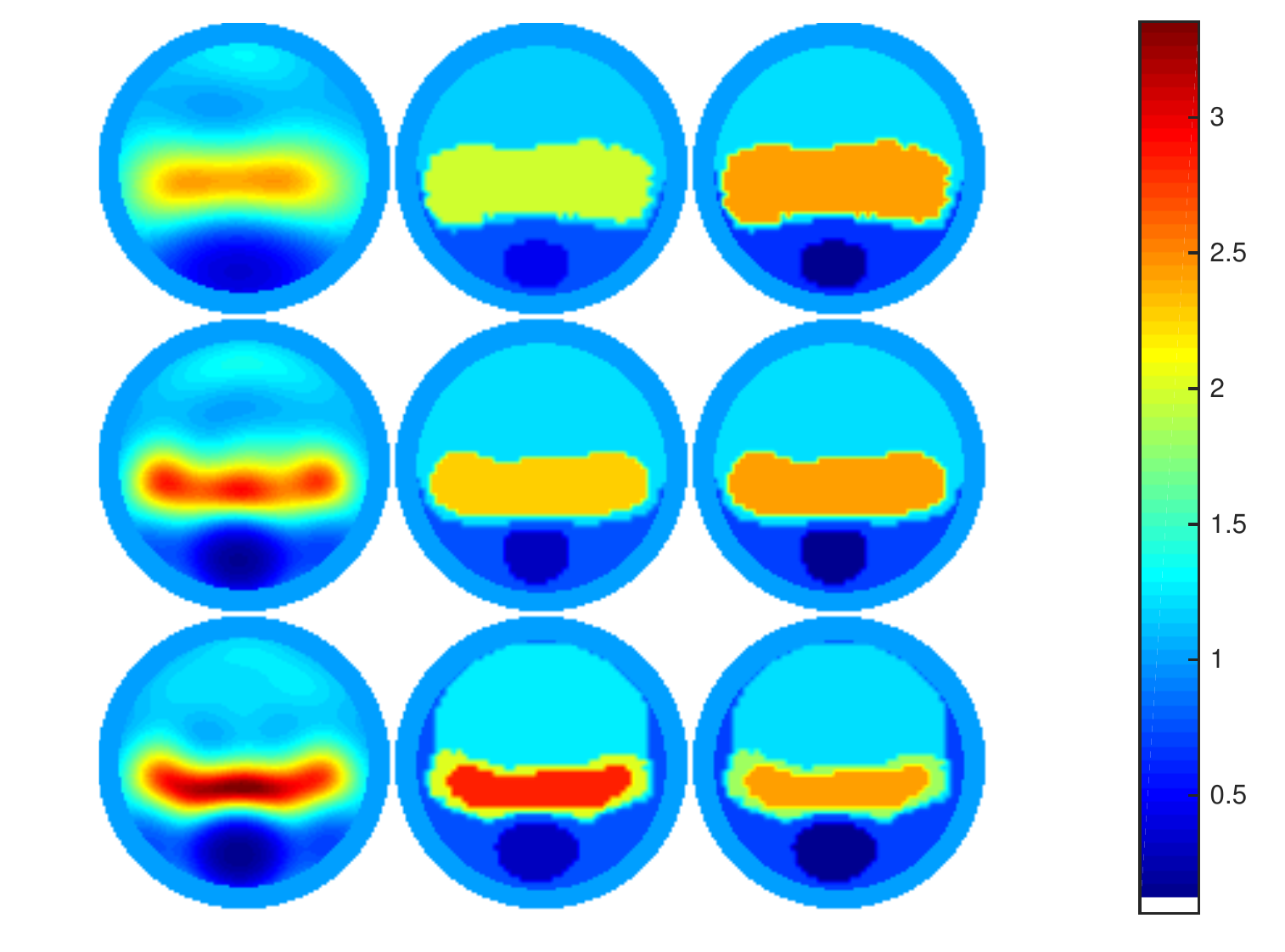}}
\put(117,40){\includegraphics[width=110pt]{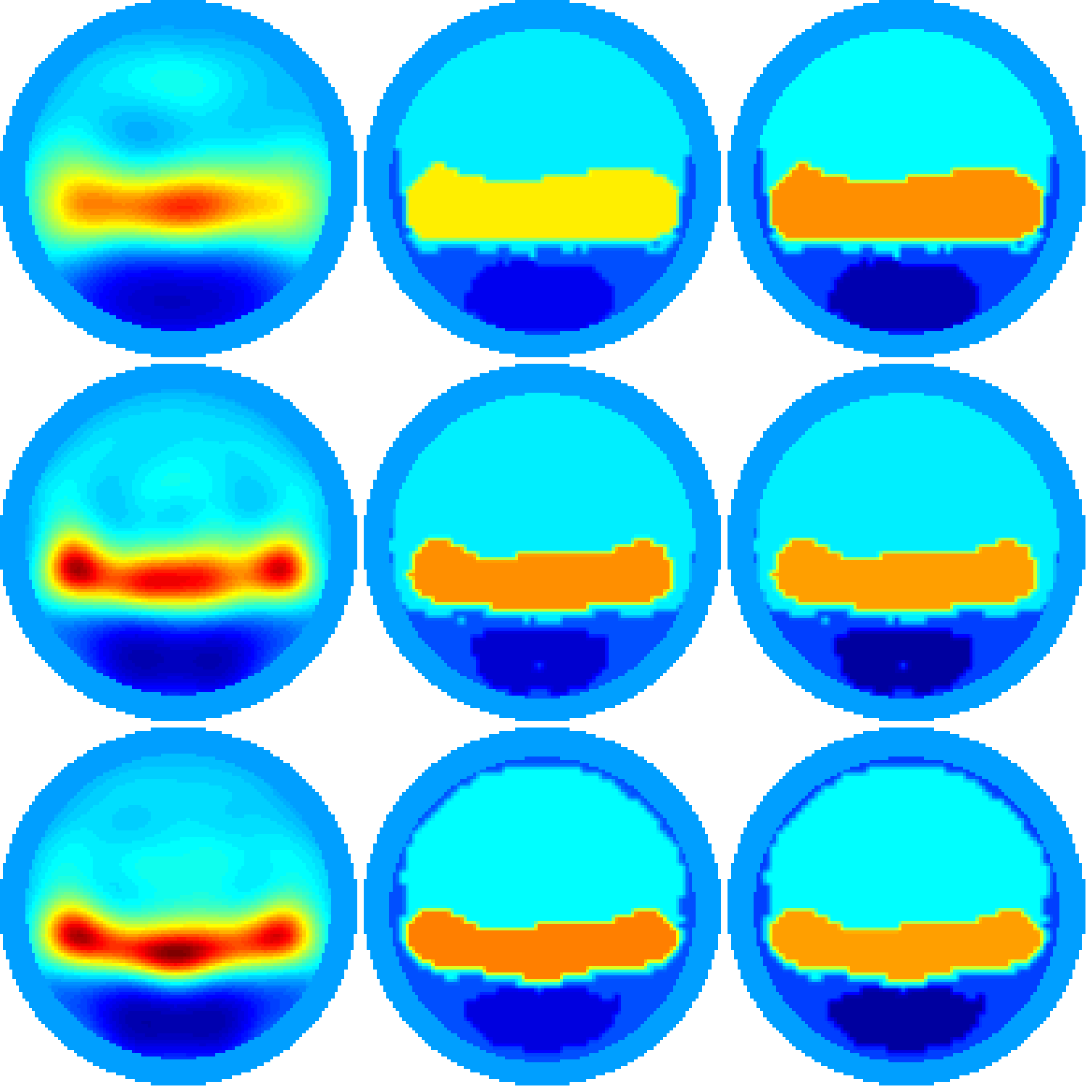}}

\put(-40,130){Iter 1}
\put(-40,95){Iter 2}
\put(-40,58){Iter 3}
\put(6,160){$\sigma_{\tiny \mbox{DB}}$}
\put(41,160){$\sigma_{\tiny \mbox{TV}}$}
\put(80,160){$\sigma_{\tiny \mbox{CE}}$}

\put(125,160){$\sigma_{\tiny \mbox{DB}}$}
\put(160,160){$\sigma_{\tiny \mbox{TV}}$}
\put(200,160){$\sigma_{\tiny \mbox{CE}}$}

\put(245,160){$\sigma_{\tiny \mbox{DB}}$}
\put(281,160){$\sigma_{\tiny \mbox{TV}}$}
\put(320,160){$\sigma_{\tiny \mbox{CE}}$}

\put(0,5){a) no added noise}
\put(130,5){b) $0.1\%$ noise}
\put(250,5){c) $0.75\%$ noise}
\end{picture}
\caption{\label{fig:all3recon_pipeline}This figure shows the real parts of the numerical approximations obtained for the industrial pipe phantom of Example~2, i.e. $\sigma_2$.  The picture consists of three parts a), b), c) corresponding to our three cases of added simulated noise in the EIT voltage data: zero added noise, noise of relative amplitude $0.1\%$ and noise of relative amplitude $0.75\%$, respectively. For the sake of comparison, the true conductivity $\sigma_2$ is displayed above the reconstructions (all on the same color scale).  For each noise level, the D-bar reconstruction $\sigma^{(j)}_{\mbox{\tiny DB}}$ (left columns), the TV sharpened image $\sigma^{(j)}_{\mbox{\tiny TV}}$ (middle columns), and the contrast adjusted TV sharpened images $\sigma^{(j)}_{\mbox{\tiny CE}}$. The first, second, and third rows correspond to the first, second and third iterations, respectively.}
\end{figure}

The reconstructed conductivities for each stage of the algorithm are displayed in Figure~\ref{fig:all3recon_pipeline}.  Note that the reconstructions are displayed on the same color scale as the original conductivity shown in Figure~\ref{fig:trueconduc} (right) for ease of comparison.   Furthermore, note that the ring of constant conductivity along the boundary (representing the thickness of the pipe) has been enforced in the reconstructions as this can also be considered \textit{apriori} information for this application.

Tables~\ref{table:Errors_0o0noise_PIPE}, \ref{table:Errors_0o1noise_PIPE}, and \ref{table:Errors_0o75noise_PIPE} show the relative $L^2$ errors and the Structural SiMilarity Index (SSIM) values for the industrial pipe phantom $\sigma_2$ of Example~2 for zero added noise, $0.1\%$ added noise, and $0.75\%$ added noise, respectively.  The error values are presented for each step of the proposed algorithm: the D-bar reconstruction $\sigma^{(j)}_{\mbox{\tiny DB}}$, the sharpened reconstruction $\sigma^{(j)}_{\mbox{\tiny TV}}$, and the contrast adjusted sharpened reconstruction $\sigma^{(j)}_{\mbox{\tiny CE}}$.

\begin{table}[h!]
\centering
{\footnotesize
\begin{tabular}{l|lllll|lll}

  \multicolumn{4}{c}{$L^2$ Relative Error} & $\quad$ & \multicolumn{4}{c}{SSIM} \\
  \multicolumn{9}{c}{}\\

  $j$ & $\sigma^{(j)}_{\mbox{\tiny DB}}$ & $\sigma^{(j)}_{\mbox{\tiny TV}}$ & $\sigma^{(j)}_{\mbox{\tiny CE}}$
   & & $j$ & $\sigma^{(j)}_{\mbox{\tiny DB}}$ & $\sigma^{(j)}_{\mbox{\tiny TV}}$ & $\sigma^{(j)}_{\mbox{\tiny CE}}$ \\

  \cline{1-4} \cline{6-9}

1 & 0.1926 & 0.1926  & 0.2157 & &
1 & 0.7097 & 0.7097 & 0.7130\\

2 & 0.1813 & 0.1985 & 0.2286 & &
2 & 0.7127 & 0.6976 & 0.6972 \\

 3 & 0.1862 & 0.2069 & 0.2228 & &
3 & 0.7054& 0.6130 & 0.6128 \\

\end{tabular}
}
\caption{\label{table:Errors_0o0noise_PIPE}The $L^2$ relative errors as well as the SSIM values for the zero added noise case for Example 2 with the pipeline phantom $\sigma_2$.}
\end{table}

\begin{table}[h!]
\centering
{\footnotesize
\begin{tabular}{l|lllll|lll}

  \multicolumn{4}{c}{$L^2$ Relative Error} & $\quad$ & \multicolumn{4}{c}{SSIM} \\
  \multicolumn{9}{c}{}\\

  $j$ & $\sigma^{(j)}_{\mbox{\tiny DB}}$ & $\sigma^{(j)}_{\mbox{\tiny TV}}$ & $\sigma^{(j)}_{\mbox{\tiny CE}}$
   & & $j$ & $\sigma^{(j)}_{\mbox{\tiny DB}}$ & $\sigma^{(j)}_{\mbox{\tiny TV}}$ & $\sigma^{(j)}_{\mbox{\tiny CE}}$ \\

  \cline{1-4} \cline{6-9}

1 & 0.1970 & 0.2176  & 0.2479 & &
1 & 0.7292 & 0.6834 & 0.6829\\

2 & 0.2279 & 0.2309 & 0.2305 & &
2 & 0.6982 & 0.6826 & 0.6837 \\

 3 & 0.2329 & 0.2404 & 0.2392 & &
3 & 0.6997& 0.5744 & 0.5653 \\

\end{tabular}
}
\caption{\label{table:Errors_0o1noise_PIPE}The $L^2$ relative errors as well as the SSIM values for the $0.1\%$ added noise case for Example 2 with the pipeline phantom $\sigma_2$.}
\end{table}

\begin{table}[h!]
\centering
{\footnotesize
\begin{tabular}{l|lllll|lll}

  \multicolumn{4}{c}{$L^2$ Relative Error} & $\quad$ & \multicolumn{4}{c}{SSIM} \\
  \multicolumn{9}{c}{}\\

  $j$ & $\sigma^{(j)}_{\mbox{\tiny DB}}$ & $\sigma^{(j)}_{\mbox{\tiny TV}}$ & $\sigma^{(j)}_{\mbox{\tiny CE}}$
   & & $j$ & $\sigma^{(j)}_{\mbox{\tiny DB}}$ & $\sigma^{(j)}_{\mbox{\tiny TV}}$ & $\sigma^{(j)}_{\mbox{\tiny CE}}$ \\

  \cline{1-4} \cline{6-9}

1 & 0.2043 & 0.2181  & 0.2617 & &
1 & 0.7123 & 0.6908 & 0.6958 \\

2 & 0.2255 & 0.2273 & 0.2434 & &
2 & 0.6873 & 0.6972 & 0.6960 \\

 3 & 0.2604 & 0.2758 & 0.2622 & &
3 & 0.6933& 0.5868 & 0.5690 \\

\end{tabular}
}
\caption{\label{table:Errors_0o75noise_PIPE}The $L^2$ relative errors as well as the SSIM values for the $0.75\%$ added noise case for Example 2 with the pipeline phantom $\sigma_2$.}
\end{table}

In the zero added noise case we see that the $L^2$ relative error in the D-bar images decreases with each iteration and the SSIM stays approximately the same.  Interestingly, the sharpened and contrast adjusted images appear to perform slightly worse under the metrics.  However, visually the sharpened images much better reflect the physical scenario (oil, water, sand) then their smooth D-bar counterparts as they contain nice clean divisions between the layers.  As the noise level increases, the method is still able to clearly distinguish between oil, water and sand at even the first iteration and thus the algorithm could be stopped there (i.e. at $\sigma_{\tiny \mbox{TV}}^{(1)}$) for each noise level.

Regarding the contrast enhanced images, recall that the approximate upper and lower bounds $C=2.5$ and $c=0.1$, respectively, were used.  In practice, closer approximations may be known (in particular for the case of oil, water, and sand) therefore improving the reconstructed values.  An investigation into the optimal parameters $c$ and $C$ and minimization scheme for the contrast adjustment is outside the scope of this introductory paper.

\section{Conclusions}\label{sec:concl}

EIT data contains information about the conductivity in an indirect, nonlinear and unstable way. Theoretically, \cite{Astala2006a} shows that infinite-precision data (the DN map) contains enough information to uniquely determine the conductivity. However, a practical data matrix $\Lambda_\sigma^\delta$ is a noisy and finite-dimensional approximation of infinite-dimensional data $\Lambda_\sigma$, and in most cases does not actually correspond to any conductivity (as $\Lambda_\sigma^\delta$ is not in the range of the forward map $\sigma\mapsto\Lambda_\sigma$). Therefore, EIT reconstruction methods need to be regularized to yield noise-robust results. Regularization is based on complementing the insufficient measurement data by {\em a priori} information about the conductivity.

Currently there are not many regularized reconstruction methods for EIT. The theory of Tikhonov regularization and related variational methods applies to a wide class of nonlinear forward maps \cite{Kaltenbacher2008}, but, alas, not to the extremely nonlinear case of EIT. For partial results, see \cite{Lechleiter2008,Jin2012a,Jin2012b}. The enclosure method for detecting convex hulls of inclusions admits a regularization analysis \cite{Ikehata2004}, but it only yields partial information (e.g. information about the locations of inclusions rather than their conductivity values). The D-bar method \cite{Knudsen2009} is the only regularized reconstruction method that produces actual conductivity images, but the reconstructions are always smooth because of a nonlinear low-pass filter involved. Indeed, the assumptions of the regularized D-bar method include continuous differentiability of the conductivity and thus the smoothing is not unexpected.

In many applications of EIT, such as nondestructive testing, the conductivity distribution can be assumed to be piecewise constant. This is approximately the case in medical imaging as well. Therefore, it is desirable to design a regularized reconstruction method producing piecewise constant images.

In this paper, a noise-robust EIT reconstruction method that always results in a piecewise constant image was both presented and tested on simulated noisy EIT data.  Therefore, this paper demonstrates that one can achieve the above goal, at least partially.  The authors note that this is an initial feasibility study only, and do not prove that the new combined method is itself a regularization strategy. However, there may be hope to prove that the reconstruction approaches the true piece-wise constant conductivity along a stable path as the data error tends to zero. Namely, the low frequencies of the reconstruction are provided by the regularized D-bar method, while the large frequencies are built based on the piecewise constant assumption. The limit frequencies between the different treatments can be assumed to tend to infinity as the noise level vanishes.

Although the scattering data in the extended annuli for $\tau^{(1)}$ and $\tau^{(2)}$ in each example is not identical to that of $\tau_{\tiny \mbox{B}}$, nonetheless, the subsequent corresponding conductivity reconstructions show marked improvements in the locations, sharpness of edges, and conductivity values of the inclusions.  Take particular note of the $0.75\%$ noise case of Example~1 where in iterations 1 and 2 the images $\sigma_{\tiny \mbox{DB}}$, $\sigma_{\tiny \mbox{TV}}$, and $\sigma_{\tiny \mbox{CE}}$ all contain a strong artifact above the heart which is not present after an additional iteration (i.e. in iteration 3).  In Example~2, we see that after a single sharpening of the original D-bar image the corresponding $\sigma_{\tiny \mbox{TV}}^{(1)}$ is enough for the goal of distinguishing between oil, water, and sand in the pipeline.  These two cases (shown for various levels of noise) demonstrate the flexibility of the algorithm: for some cases it is appropriate to be applied iteratively for improvements and to remove artifacts, whereas in other cases a single sharpening iteration is adequate.


Additional modifications to the proposed method can be easily applied.  In particular, in lieu of fixing the maximal number of iterations $J$ beforehand, alternative stopping criteria could be applied.  E.g., for $j>1$ instead of asking above if $j=J$, return $\sigma_{\mbox{\tiny CE}}^{(j)}$ as the final image if
\[\norm{\sigma_{\mbox{\tiny CE}}^{(j)}-\sigma_{\mbox{\tiny CE}}^{(j-1)}}_{l^2} / \norm{\sigma_{\mbox{\tiny CE}}^{(j)}}_{l^2} <thresh\]
for some specified threshold $thresh$.  Another option concerns the radii $R$ and $\widetilde{R}$.  In the examples presented here, the radii $R$ and $\widetilde{R}$ are fixed throughout the algorithm for each data set $\Lambda_\sigma$. An alternative approach could compute
the Beltrami scattering data on progressively larger annuli so
that in Step~(e) of the algorithm in Figure~\ref{fig:flowchart}, the new scattering data
$\widetilde{\tau}^{(j)}(k)$ is computed for $R-1<|k|\leq
\widetilde{R} + (j-1)\Delta R$ given a fixed stepsize $\Delta R>0$
($j\geq 1$). Such approaches, while interesting are outside the scope of this work.

\section*{Acknowledgments}

\noindent
S.~J. Hamilton, J.~M. Reyes, and S. Siltanen were supported by the Academy of Finland (Finnish Centre of Excellence in Inverse Problems Research 2012–2017,
decision number 250215).  S.~J. Hamilton was additionally supported by SalWe Research Program for Mind and Body (Tekes - the Finnish Funding Agency for Technology
and Innovation grant 1104/10).
J.~M. Reyes was additionally supported by the Engineering and Physical
Sciences Research Council (EPSRC), reference EP/K024078/1.  X. Zhang was supported by NSFC11101277 and  NSFC91330102.


\bibliographystyle{plain}

\bibliography{bibliographyRefs}

\end{document}